\documentclass[aps,pre,10pt,showpacs,floatfix,onecolumn,nofootinbib,a4paper]{revtex4-2}

\usepackage{amssymb, amsmath, amsthm, mathtools}
\usepackage{mathrsfs}
\usepackage{hyperref,bookmark}
\usepackage{graphicx}
\usepackage{epsfig,color}
\usepackage{float}
\usepackage{caption}
\usepackage{subcaption}
\usepackage{verbatim}
\usepackage{dcolumn}
\usepackage{tikz}
\usepackage{url}
\usepackage{cases}
\usepackage[utf8x]{inputenc}
\usepackage[english]{babel}
\usepackage{epsfig}
\usepackage{listings}
\usepackage{paralist}
\usepackage{mdwlist}
\usepackage{bm}
\usepackage{dsfont}
\usepackage{oldgerm}
\usepackage{geometry}
\usepackage{relsize}
\usepackage{afterpage}
\usepackage{lmodern}
\usepackage{tabularx}
\usepackage{cancel}
\usepackage{dutchcal} 
\usepackage{natbib}

\newcommand{\e}{\bm{e}}

\newcommand{\x}{\bm{x}}
\newcommand{\y}{\bm{y}}

\newcommand{\I}{\mathbf{I}}

\newcommand{\sphere}{\mathbb{S}^2}

\newcommand{\curvature}{(\nablas\normal)}

\newcommand{\nablas}{\nabla\!_\mathrm{s}}

\newcommand{\tr}{\operatorname{tr}}

\newcommand{\trans}{^\mathsf{T}}

\newcommand{\surface}{\mathscr{S}}
\newcommand{\normal}{\bm{\nu}}
\newcommand{\dd}{\mathrm{d}}

\newcommand{\euclid}{\mathscr{E}}

\newcommand{\transl}{\mathscr{V}}
\newcommand{\tplane}{\mathscr{T}}

\newcommand{\proj}{\mathbf{P}(\normal)}
\newcommand{\Proj}{\mathbf{P}}
\newcommand{\A}{\mathbf{A}}

\renewcommand{\H}{\mathbf{H}}

\newcommand{\U}{\mathbf{U}}
\newcommand{\V}{\mathbf{V}}
\newcommand{\W}{\mathbf{W}}
\newcommand{\R}{\mathbf{R}}
\newcommand{\av}{\bm{a}}

\newcommand{\cv}{\bm{c}}
\newcommand{\dv}{\bm{d}}
\newcommand{\rv}{\bm{r}}

\newcommand{\wn}{\lvert w\rvert}
\newcommand{\wsn}{\lvert w^\ast\rvert}
\newcommand{\hn}{\lvert h\vert}
\newcommand{\hpn}{\lvert h'\rvert}

\newcommand{\ws}{W_\mathrm{s}}
\newcommand{\wdr}{W_\mathrm{d}}
\newcommand{\wb}{W_\mathrm{b}}

\newcommand{\nablast}{\nabla^\ast\!\!\!_\mathrm{s}}

\newcommand{\framee}{(\e_1,\e_2,\e_3)}
\newcommand{\ex}{\mathrm{e}}

\begin{document}
\latintext
\title{Neutral Deformation Modes of Minimal Surfaces}
\author{Andr\'e M. Sonnet}
\email{andre.sonnet@strath.ac.uk}
\affiliation{Department of Mathematics and Statistics, 26 Richmond Street, Glasgow, G1 1XH, U.K}
\author{Epifanio G. Virga}
\email{eg.virga@unipv.it}
\affiliation{Dipartimento di Matematica, Universit\`a di Pavia, Via Ferrata 5, 27100 Pavia, Italy}

\date{\today}

\begin{abstract}
	Stretching, drilling, and bending are the independent deformation modes of a thin shell, each of which has an individual energy content. When the energy content of a mode vanishes, that mode is \emph{neutral}. We characterize all neutral modes of deformation of minimal surfaces into minimal surfaces. A hierarchy is found among these: a stretching neutral mode (which is an isometry) is also drilling neutral, and a drilling neutral mode is also bending neutral. Thus, all isometries of a minimal surface are globally neutral and give rise to \emph{soft elasticity}. More generally, all minimal surfaces can be classified relative to a reference one in terms of three energy  contents, which can be given in closed form.
\end{abstract}

\maketitle

\section{Introduction}\label{sec:intro}
\emph{Soft} thin shells are two-dimensional material bodies, likely to be more responsive than others to \emph{drilling} deformations, which are produced by local rotations pivoting about the local normal.

In the variational theory that we  proposed  in \cite{sonnet:variational}, one such shell is represented as a material surface $\surface$. The strain energy density features \emph{three} distinct contributions: together with the customary \emph{stretching} energy, formulated in terms of the metric changes induced on $\surface$, and the \emph{bending} energy, related to the changes in curvature, a novel \emph{drilling} energy is introduced, expressed in terms of changes in the \emph{spin connector} of an appropriate moving frame.

Our theory for soft thin shells is formally \emph{direct}, as it is not \emph{derived} from three-dimensional elasticity,\footnote{We refer the reader to Sects.\,212,\,213 of \cite{truesdell:classical} for a classical, but still effective summary of direct and derived approaches to the elastic theory of shells. Several applications of the direct approach can also be found in the books \cite{antman:nonlinear,rubin:cosserat}. Finally, \cite{tomassetti:coordinate} is a more recent tutorial on direct theories advocating the same coordinate-free formalism adopted here.} yet it is still motivated by an average process; but neither are  averages taken across the thickness of a three-dimensional body, nor are three-dimensional strains being averaged about.

The realms of direct and derived theories for shells were connected by the theory of the Cosserat brothers \cite{cosserat:theorie,cosserat:theorie_livre}, who first used \emph{directors} (see also \cite{ericksen:exact} and \cite{cohen:nonlinear_directed}) to recover for lower-dimensional bodies the information lost in the neglected dimensions, as shown in the classical essay \cite{naghdi:theory}.

The avenue pursued in \cite{sonnet:variational} is different. While two-dimensional stretching measures are classically derived from the theory of membranes, drilling and bending measures require motivation; this comes from averaging relative local rotations on a moving tangent disc. Thus, averaging is tangential instead of transversal, and it only involves rotations instead of full deformations gradients.

This approach led us in \cite{sonnet:variational} to identify the strain energy density of soft thin shells as the sum of three independent \emph{pure}  measures of stretching, drilling, and bending, with each individual measure vanishing on all other companion modes.\footnote{A deformation measure is \emph{pure} if it is affected by a single mode; see, for example, the bending measures studied in \cite{acharya:nonlinear,wood:contrasting,vitral:dilation,acharya:mid-surface}.}

As common to other theories, including director theories whose single director can be identified with the unit normal to $\surface$, see for example \cite{cohen:nonlinear,hilgers:Gauss,sanders:nonlinear,silhavy:direct}, the proposed strain energy density depends on both the first and second deformation gradients, and so our averaging process produces effectively a second-grade theory, although it does so only through the gradient of the \emph{polar} rotation.\footnote{This is the rotation in three-dimensional space extracted from the deformation gradient through its polar decomposition.}

In a broad sense, the theory for soft thin shells in \cite{sonnet:variational} appears as the natural extension of Antman's theory for elastic rods \cite{antman:general}, where the splitting in pure deformation modes made perhaps its first appearance.

Any new theory is likely (and also lucky) to have supporters and detractors (most go simply unnoticed). Especially on crowded grounds, as is that of shell theories, one needs direct experience to move forward. It is the purpose of this paper to put our theory to the test. Before attempting to solve specific boundary value problems with it, we found it expedient to see how the energy splitting characteristic of its construction behaves on a wealth of shell shapes and their deformations; of course, the wider the better. We chose the whole class of \emph{minimal surfaces}  to describe the three deformation modes and the corresponding energies. The multitude of shapes in this class sufficed to outline a hierarchical pattern, which confirms a rich underlying energy landscape, including instances of \emph{soft elasticity} (with a continuum of ground states). By considering all possible deformations of a minimal surface, we complete the programme started in \cite{sonnet:bending-neutral}, where only bending neutral deformations of minimal surfaces were considered.

The rest of the paper is organized as follows. In Sect.~\ref{sec:minimal}, we recall the classical representation of minimal surfaces in terms of the Weierstrass function and show how a deformation of a minimal surface into another can be described by a conformal mapping on the complex plane. Section~\ref{sec:modes} is concerned with the identification of the different independent modes and the representation of their energy contents. In Sect.~\ref{sec:neutral}, we isolate deformations of a minimal surface for which a single energy content vanishes selectively; these are the \emph{neutral} modes, which can be identified in great generality, as all energy contents are given in closed form. Explicit examples are also provided in which both Bour's and Enneper's surfaces undergo selective deformations. Finally, in Sect.~\ref{sec:conclusion}, we summarize the conclusions of this study. Two appendices containing a number of mathematical details close the paper.

\section{Minimal Surfaces}\label{sec:minimal}

\subsection{Weierstrass representation}\label{sec:Weierstrass}
Let $w:=u+iv\in\mathbb{C}$ be a complex representation of the $(u,v)$ plane. A minimal surface $\surface$ such that the normal $\normal$ regarded as a map from a domain $\Omega$ in the plane into the unit sphere $\sphere$ is injective, but does not cover the whole of $\sphere$, can be represented (to within a translation) by a \emph{holomorphic} function $F:\Omega\to\mathbb{C}$ as
\begin{equation}
\label{eq:Weierstrass_representation}
\rv(w)=\Re\left(\frac12\int(1-w^2)F(w)\dd w\,\e_1+\frac{i}{2}\int(1+w^2)F(w)\dd w\,\e_2+\int wF(w)\dd w\,\e_3\right),
\end{equation}
where $\framee$ is a Cartesian frame in $\transl$ and $\Re$ denotes the real part (see, for example, \cite[p.\,117]{dierkes:minimal}).\footnote{Here, $\transl$ is the translation space associated with three-dimensional Euclidean space $\euclid$, see, for example, \cite[p.\,324]{truesdell:first}, where these basic geometric structures are further illuminated.} $F$ is also called the \emph{Weierstrass function} that represents $\surface$ (see also \cite{weierstrass:untersuchungen_berlin,weierstrass:untersuchungen_cambridge}).

A number of interesting consequences follow from \eqref{eq:Weierstrass_representation}, which can easily be interpreted geometrically.

We first write $F$ as
\begin{equation}
	\label{eq:F_representation}
	F=\ex^{\Phi+i\chi}=\ex^\Phi(\cos\chi+i\sin\chi).
\end{equation}
Since $F$ is holomorphic, the Cauchy-Riemann relations imply that both $\Phi$ and $\chi$ are harmonic in $\Omega$. Explicit computations give%
\footnote{Here and in the following a comma denotes partial differentiation with respect to the indicated variable.}
\begin{subequations}\label{eq:Weierstrass_metric_vectors}
\begin{align}
\rv_{,u}&=\ex^\Phi\Big\{\Big[uv\sin\chi+\frac12(1-u^2+v^2)\cos\chi\Big]\e_1-\Big[uv\cos\chi+\frac12(1+u^2-v^2)\sin\chi\Big]\e_2\nonumber\\&+(u\cos\chi-v\sin\chi)\e_3\Big\},\label{eq:Weierstrass_metric_vector_u}\\
\rv_{,v}&=\ex^\Phi\Big\{\Big[uv\cos\chi-\frac12(1-u^2+v^2)\sin\chi\Big]\e_1+\Big[uv\sin\chi-\frac12(1+u^2-v^2)\cos\chi\Big]\e_2\nonumber\\&-(u\sin\chi+v\cos\chi)\e_3\Big\},\label{eq:Weierstrass_metric_vector_v}
\end{align}
which satisfy
\begin{equation}\label{eq:Weierstrass_metric_vector_u_v}
\lvert\rv_{,u}\rvert=\lvert\rv_{,v}\rvert=\frac12\ex^\Phi(1+u^2+v^2)=\frac12\ex^\Phi(1+\wn^2)\quad\text{and}\quad\rv_{,u}\cdot\rv_{,v}=0.
\end{equation}
\end{subequations}
Let $\e_u$ and $\e_v$ denote the unit vectors of $\rv_{,u}$ and $\rv_{,v}$, we have that
\begin{equation}\label{eq:Weierstrass_frame}
\e_u=\frac{\rv_{,u}}{\lvert\rv_{,u}\rvert},\quad
\e_v=\frac{\rv_{,v}}{\lvert\rv_{,v}\rvert}.
\end{equation}
They form an orthonormal frame $(\e_u,\e_v)$ on the tangent plane $\tplane_{\rv(w)}$ on $\surface$; the unit normal $\normal$ is then obtained from
\begin{equation}\label{eq:Weierstrass_normal}
\normal=\e_u\times\e_v=\frac{1}{1+\wn^2}\{2u\e_1+2v\e_2+(\wn^2-1)\e_3\}.
\end{equation}

Consider now a trajectory $(u(t),v(t))$ in $\Omega$ parametrised by $t$; $\rv$  maps it into a trajectory on $\surface$ for which
\begin{equation}
\label{eq:Weierstrass_r_dot}
\dot{\rv}=\lvert\rv_{,u}\rvert(\dot{u}\e_u+\dot{v}\e_v),
\end{equation}
where a superimposed dot denotes differentiation with respect  to $t$. Correspondingly,
\begin{equation}
	\label{eq:Weierstrass_normal_dot}
	\dot{\normal}=\dot{u}\normal_{,u}+\dot{v}\normal_{,v}=\curvature\dot{\rv},
\end{equation}
where $\curvature$ is the \emph{curvature tensor} of $\surface$.
By differentiating with respect to $u$ and $v$ both sides of \eqref{eq:Weierstrass_normal} and projecting both $\normal_{,u}$ and $\normal_{,v}$ in the frame $(\e_u,\e_v)$, we obtain that 
\begin{subequations}\label{eq:Weierstrass_partial_normal}
\begin{align}
\normal_{,u}&=\frac{2}{1+\wn^2}(\cos\chi\e_u-\sin\chi\e_v),	\label{eq:Weierstrass_partial_normal_u}\\
\normal_{,v}&=-\frac{2}{1+\wn^2}(\sin\chi\e_u+\cos\chi\e_v).\label{eq:Weierstrass_partial_normal_v}
\end{align}
\end{subequations}
Inserting \eqref{eq:Weierstrass_partial_normal} and \eqref{eq:Weierstrass_r_dot} into \eqref{eq:Weierstrass_normal_dot}, with the aid of \eqref{eq:Weierstrass_metric_vector_u_v}, since $\dot{u}$ and $\dot{v}$ are arbitrary, we arrive at the following representation of the curvature tensor in terms of the functions $\Phi$ and $\chi$,
\begin{equation}\label{eq:Weierstrass_curvature_tensor}
\nablas\normal=\frac{4\ex^{-\Phi}}{(1+\wn^2)^2}
\left[\cos\chi(\e_u\otimes\e_u-\e_v\otimes\e_v)-
\sin\chi(\e_u\otimes\e_v+\e_v\otimes\e_u)
\right].
\end{equation}
Computing the trace and determinant of the curvature tensor in this form is straightforward and shows that the mean curvature is indeed zero, $2H=\tr\curvature=0$, and that the Gaussian curvature is
\begin{equation}\label{eq:Weierstrass_curvature}
K=\det\curvature=-\frac{16\ex^{-2\Phi}}{(1+\wn^2)^4}.
\end{equation}

\subsection{Surface Kinematics}\label{sec:kinematics}
\begin{figure}[h]
\centering\includegraphics[width=.8\textwidth]{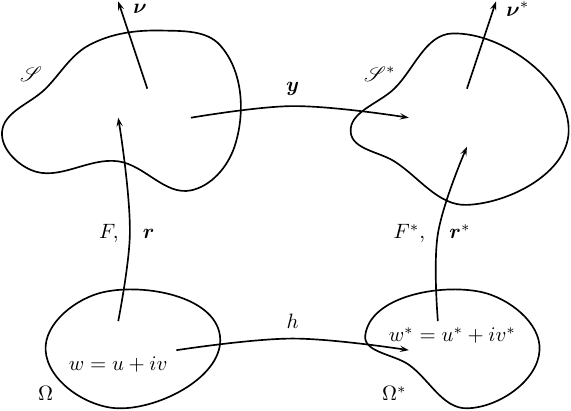}
\caption{Let two minimal surfaces $\surface$ and $\surface^\ast$ be represented by Weierstrass functions $F$ and $F^\ast$ with simply connected open domains $\Omega$ and $\Omega^\ast$. Let $h$ be a conformal mapping such that $\Omega^\ast=h(\Omega)$.
We consider the deformation $\y$ that takes $\surface$ into the surface $\surface^{\ast}=\y(\surface)$ and satisfies $\y\circ\rv =\rv^\ast\circ h$ where $\rv$ and $\rv^\ast$
are the mappings produced by $F$ and $F^\ast$ respectively according to \eqref{eq:Weierstrass_representation}.}
\label{fig:mappings}
\end{figure}

Consider locally two arbitrary minimal surfaces $\surface$ and $\surface^\ast$, represented by Weierstrass functions
$F=\ex^{\Phi+i\chi}$ and $F^\ast=\ex^{\Phi^\ast+i\chi^\ast}$ defined on domains $\Omega$ and $\Omega^\ast$, respectively.
We assume that there exists a conformal mapping $h$ such that%
\footnote{So $h$ is analytic and $\hpn\ne0$.}
\begin{equation}\label{eq:mappingh}
\Omega^\ast=h(\Omega). 
\end{equation}
If $\Omega=\Omega^\ast$ as in \cite{sonnet:bending-neutral} the identity will work, and this includes the special case $\Omega=\mathbb{C}$. More generally, if
both domains are open, simply connected, and strictly smaller than $\mathbb{C}$,
Riemann's mapping theorem guarantees the existence of a biholomorphic
conformal mapping between them (see, for example \cite[p. 283]{rudin:analysis}).

Let $\rv$ and $\rv^\ast$ be the mappings produced by the two Weierstrass functions $F$ and $F^\ast$ according to \eqref{eq:Weierstrass_representation}.
We are interested in the deformation $\y:\surface\rightarrow\surface^\ast$ that satisfies
\begin{equation}\label{eq:deformationMapping}
\y\circ\rv =\rv^\ast\circ h,
\end{equation}
the situation depicted in Figure~\ref{fig:mappings}.
In other words, with $w=u+iv\in\Omega$ and letting
\begin{equation}
	\label{eq:h_defintion}
	w^\ast:=h(w)=u^\ast+iv^\ast\in\Omega^\ast
\end{equation}
we have
\begin{equation}\label{eq:deformation}
\y(\rv(u+iv))=\rv^\ast(h(u+iv)).
\end{equation}

To simplify notation, let us call
\begin{equation}\label{eq:u_v_star}
u^\ast=h_u(u,v)\quad\text{and}\quad v^\ast =h_v(u,v).
\end{equation}
The Cauchy-Riemann equations for $h$ then read as
\begin{equation}\label{eq:CRequations}
 h_{u,u}=h_{v,v}\quad\text{and}\quad h_{u,v}=-h_{v,u}.
\end{equation}
We note for later use that
\begin{equation}\label{eq:hprime}
h'=h_{u,u}+i h_{v,u}\quad\text{and so}\quad
\lvert h'\rvert^2=h_{u,u}^2+ h_{v,u}^2.
\end{equation}

We introduce two orthonormal frames on $\surface^\ast$. One
is provided naturally by the starred version of \eqref{eq:Weierstrass_frame},
\begin{equation}\label{eq:DoubleStarredBasis}
\e^\ast_{u^\ast}=\frac{\rv_{,u^\ast}^\ast}{\lvert\rv_{,u^\ast}^\ast\rvert},\quad
\e^\ast_{v^\ast}=\frac{\rv_{,v^\ast}^\ast}{\lvert\rv_{,v^\ast}^\ast\rvert}
\end{equation}
with
\begin{equation}
 \lvert\rv^\ast_{,u^\ast}\rvert=\lvert\rv^\ast_{,v^\ast}\rvert=\frac{1}{2}e^{\Phi^\ast}(1+\wsn^2)=\frac{1}{2}e^{\Phi^\ast}(1+\hn^2).
\end{equation}

However, in order to analyse the deformation $\y$ of $\surface$ into $\surface^\ast$, we want to parametrise also the starred surface in terms of the unstarred parameters $u$ and $v$.
By taking the derivatives of \eqref{eq:deformation} with respect to $u$ and $v$ and using the chain
rule on both sides, we obtain
\begin{subequations}\label{eq:mapr}
\begin{align}\label{eq:mapru}
 (\nablas\y)\rv_{,u}&=\rv^\ast_{,u}=\rv^\ast_{,u^\ast}h_{u,u}+\rv^\ast_{,v^\ast}h_{v,u}\\
 \label{eq:maprv}
(\nablas\y)\rv_{,v}&=\rv^\ast_{,v}=\rv^\ast_{,u^\ast}h_{u,v}+\rv^\ast_{,v^\ast}h_{v,v},
\end{align}
\end{subequations}
where $\nablas\y$ is the surface deformation gradient of $\y$.

Computing the norms of the right-hand sides of \eqref{eq:mapr} shows that
\begin{equation}\label{eq:rStarunorm}
\lvert\rv^\ast_{,u}\rvert=\lvert\rv^\ast_{,v}\rvert=\sqrt{\lvert\rv^\ast_{,u^\ast}\rvert^2 h_{u,u}^2
+\lvert\rv^\ast_{,v^\ast}\rvert^2 h_{v,u}^2}=
 \lvert\rv^\ast_{,u^\ast}\rvert\,\lvert h'\rvert,
\end{equation}
where we have used the Cauchy-Riemann equations for $h$ along with \eqref{eq:hprime}.
We then obtain the desired second orthonormal frame
$(\e^\ast_u, \e^\ast_v$) of unit vectors on the tangent plane to $\surface^\ast$
by defining
\begin{equation}\label{eq:StarredBasis}
\e_u^\ast=\frac{\rv_{,u}^\ast}{\lvert\rv_{,u}^\ast\rvert},\quad
\e_v^\ast=\frac{\rv_{,v}^\ast}{\lvert\rv_{,v}^\ast\rvert}.
\end{equation}
Using \eqref{eq:mapr} and \eqref{eq:rStarunorm}, it can be seen that the two frames on $\surface^\ast$ are related by
\begin{subequations}\label{eq:StarredBases}
\begin{align}
\e_u^\ast&=\frac{1}{\hpn}\left(
h_{u,u}\e^\ast_{u^\ast}+h_{v,u}\e^\ast_{v^\ast}\right),\label{eq:estaru} \\
\e_v^\ast&=\frac{1}{\hpn}\left(
h_{u,v}\e^\ast_{u^\ast}+h_{v,v}\e^\ast_{v^\ast}\right),
\end{align}
\end{subequations}
which can be inverted to yield
\begin{subequations}\label{eq:StarredBasesInverted}
\begin{align}
\e_{u^\ast}^\ast&=\frac{1}{\hpn}\left(
h_{v,v}\e^\ast_{u}-h_{v,u}\e^\ast_{v}\right), \\
\e_{v^\ast}^\ast&=\frac{1}{\hpn}\left(
-h_{u,v}\e^\ast_{u}+h_{u,u}\e^\ast_{v}\right).
\end{align}
\end{subequations}

By the polar decomposition theorem for the deformation of surfaces~\cite{man:coordinate}, the deformation gradient $\nablas\y$
can be written as
\begin{equation}\label{eq:polar_decomposition}
\nablas\y=\R\U=\V\R,
\end{equation}
where $\R$ is a three-dimensional rotation and $\U$ and $\V$ are positive definite symmetric stretching tensors that operate on the tangent planes to
$\surface$ and $\surface^\ast$, respectively.

We want to find explicit expressions for $\R$, $\U$, and $\V$, which will
allow us in Section~\ref{sec:neutral} to compute the individual deformation
modes described in Section~\ref{sec:modes}. Our starting point is the fact
seen from \eqref{eq:mapr} that the
deformation gradient maps $\rv_{,u}$ to $\rv^\ast_{,u}$ and
$\rv_{,v}$ to $\rv^\ast_{,v}$.

So the deformation gradient is
\begin{equation}\label{eq:deformation_gradient}
\nablas\y=\frac{\lvert\rv_{,u}^\ast\rvert}{\lvert\rv_{,u}\rvert}(\e_u^\ast\otimes\e_u+\e_v^\ast\otimes\e_v)
=e^{\Phi^\ast-\Phi}\frac{1+\hn^2}{1+\wn^2}\hpn(\e_u^\ast\otimes\e_u+\e_v^\ast\otimes\e_v)
\end{equation}
The rotation $\R$ in \eqref{eq:polar_decomposition} takes
$(\e_u,\e_v,\normal)$ into $(\e_u^
\ast,\e_v^\ast,\normal^\ast)$, and so%
\footnote{We have $\normal=\e_u\times\e_v$ and $\normal^\ast=\e^\ast_{u^\ast}\times\e^\ast_{v^\ast}$ as in \eqref{eq:Weierstrass_normal}. From \eqref{eq:mapr} it follows that $\rv_{,u}^\ast\times\rv_{,v}^\ast=\hpn^2
\rv_{,u^\ast}^\ast\times\rv_{,v^\ast}^\ast$, and so also
$\normal^\ast=\e^\ast_u\times\e^\ast_v$, the two frames on $\surface^\ast$ share the same orientation.}
\begin{equation}\label{eq:rotationDyadic}
\R=\e_u^\ast\otimes\e_u+\e_v^\ast\otimes\e_v+\normal^\ast\otimes\normal.
\end{equation}
Thus, by inspection,
\begin{equation}\label{eq:stretchingU}
\U=\frac{\lvert\rv_{,u}^\ast\rvert}{\lvert\rv_{,u}\rvert}\proj =e^{\Phi^\ast-\Phi}\frac{1+\hn^2}{1+\wn^2}\hpn\, \proj
\end{equation}
and
\begin{equation}\label{eq:stretching}
\V=\frac{\lvert\rv_{,u}^\ast\rvert}{\lvert\rv_{,u}\rvert}\Proj(\normal^\ast) =e^{\Phi^\ast-\Phi}\frac{1+\hn^2}{1+\wn^2}\hpn\, \Proj(\normal^\ast),
\end{equation}
where
\begin{equation}\label{eq:projectors}
\proj=\I-\normal\otimes\normal\quad\text{and}\quad
\Proj(\normal^\ast)=\I-\normal^\ast\otimes\normal^\ast
\end{equation}
are the projectors onto the tangent planes of
$\surface$ and $\surface^\ast$.

A general remark should be made about the notation employed here and in the following. Starred quantities, such as $\Phi^\ast$ and $\chi^\ast$, depend on $w^\ast$; they can also be regarded as functions of $w$, if use is made of \eqref{eq:h_defintion}. When both starred and unstarred quantities appear in a formula, as in \eqref{eq:deformation_gradient} and a number of other places below, we tacitly mean that they are expressed through \eqref{eq:h_defintion} in terms of one and the same variable, which is usually $w$. This is, for example, explicitly spelled out in the following expression for the Gaussian curvature $K^\ast$ of $\surface^\ast$,
\begin{equation}
	\label{eq:K_star}
	\frac{K^\ast(h(w))}{K(w)}=e^{-2(\Phi^\ast(h(w)))-\Phi(w)}\frac{(1+|w|^2)^4}{(1+|h(w)|^2)^4},
\end{equation}
which is easily obtained from \eqref{eq:Weierstrass_curvature} and its analogue for $\surface^\ast$. To avoid typographical clutter, we shall seldom go into the detail displayed in \eqref{eq:K_star}.

\section{Modes of Deformation}\label{sec:modes}
In \cite{sonnet:variational} we identified three distinct modes in the deformation of a thin shell: a metric one pertaining to a stretching of the material in the tangent plane, a drilling one pertaining to a rotation about the local surface normal, and a bending one pertaining
to a rotation about a vector tangent to the surface.

A starting point for arriving at elastic energy measures for the three
distinct modes is the polar decomposition \eqref{eq:polar_decomposition}.
In an isometry, the stretching tensors $\U$ and $\V$ are simply
the projectors onto the tangent planes of $\surface$ and $\surface^\ast$,
and so $\ws$ given by
\begin{equation}\label{eq:w_s}
\ws=\lvert\U-\proj\rvert^2=\lvert\V-\Proj(\normal^{\ast})\rvert^2
\end{equation}
is a natural choice for a pure measure of stretching.

The rotation $\R$ in \eqref{eq:polar_decomposition} can have contributions from both drilling and bending. It was proved in \cite{sonnet:bending-neutral} that any rotation can be uniquely decomposed into one about a given unit vector $\e$ and one about a vector perpendicular to $\e$. Hence,
letting this $\e$ be the normal $\normal$, the rotation
$\R$ can be written in the form
\begin{equation}
\R=\R^b(\alpha_b)\R^d(\alpha_d)
\end{equation}
where $\R^b$ is a rotation by an angle $\alpha_b$ about a direction in the tangent plane and
$\R^d$ is a rotation by an angle $\alpha_d$ about the surface normal. One might be tempted to use the two rotation angles $\alpha_b$ and $\alpha_d$ to construct measures of bending and drilling, but unfortunately
they are not frame invariant. A uniform rotation of the shell performed after the deformation $\y$ would change the angles in an intricate way~\cite{sonnet:variational}.

A way forward was proposed in \cite{sonnet:variational}. Instead of using the angles $\alpha_b$ and $\alpha_d$ directly, local averages over relative rotations are computed in small discs on the tangent plane%
\footnote{To be precise, the averages computed in \cite{sonnet:variational}
are of $\tan^2(\alpha_b/2)$ and $\tan^2(\alpha_d/2)$.}. Using these averages, pure measures of drilling and bending can then be obtained in terms of the third-rank tensor
\begin{equation}
\label{eq:H_definition}
\H:=\R\trans\nablas\R,
\end{equation}
the skew-symmetric tensor $\W(\normal)$ associted with the normal
$\normal$ via
\begin{equation}
\W(\normal)\av=\normal\times\av
\end{equation}
for all vectors $\av$, and the curvature tensor $\nablas\normal$.
A pure measure of drilling is given by
\begin{equation}\label{eq:w_d}
\wdr:=\lvert\W(\normal)\circ\H\rvert^2
\end{equation}
and a pure measure of bending by
\begin{equation}\label{eq:w_b}
\wb:=
\left(
\lvert\H\rvert^2-\frac12 \lvert\W(\normal)\circ\H\rvert^2-4\normal\cdot\H\circ\nablas\normal
\right)^2.
\end{equation}
Here, the two vector-valued products of a third- and a second-rank tensor
indicated by $\circ$ can be represented in component form as $\A\circ \H=A_{ij}H_{ijk}\e_k$ and
$\H\circ \A=H_{ijk}A_{jk}\e_i$, where the convention of summing over repeated indices has been used.

The elastic energy density associated with a deformation of a shell would then be
\begin{equation}\label{eq:W}
W=\frac12\mu_\mathrm{s}\ws+\frac12\mu_\mathrm{d}\wdr+\frac14\mu_\mathrm{b}\wb,
\end{equation}
where the positive scalars $\mu_\mathrm{s}$, $\mu_\mathrm{d}$, and $\mu_\mathrm{b}$ are the \emph{stretching}, \emph{drilling}, and \emph{bending} moduli, respectively.

\subsection{Geometric Representation}\label{sec:representation}
Both $\wdr$ and $\wb$ can be computed explicitly in terms of systems of connectors for the moving frames $(\e_u,\e_v,\normal)$ and $(\e_u^\ast,\e_v^\ast,\normal^\ast)$ on the surfaces $\surface$ and $\surface^\ast$. The connectors are fields $(\cv,\dv_u,\dv_v)$ and $(\cv^\ast,\dv_u^\ast,\dv_v^\ast)$, tangent to $\surface$ and $\surface^\ast$, that represent the derivatives of the frames' basis vectors.
Explicitly, they obey the equations
\begin{equation}\label{eq:connectors}
  \begin{cases}
    \nablas\e_u=\e_v\otimes\cv+\normal\otimes\dv_u,\\
    \nablas\e_v=-\e_u\otimes\cv+\normal\otimes\dv_v,\\
    \nablas\normal=-\e_u\otimes\dv_u-\e_v\otimes\dv_v,
  \end{cases}
\qquad
  \begin{cases}
    \nablast\e_u^\ast=\e_v^\ast\otimes\cv^\ast+\normal^\ast\otimes\dv_u^\ast,\\
    \nablast\e_v^\ast=-\e_u^\ast\otimes\cv^\ast+\normal^\ast\otimes\dv_v^\ast,\\
    \nablast\normal^\ast=-\e_u^\ast\otimes\dv_u^\ast-\e_v^\ast\otimes\dv_v^\ast,
  \end{cases}
\end{equation}
where $\nablast$ denotes the surface gradient on $\surface^\ast$. We call $\cv$ and $\cv^\ast$ the \emph{spin} connectors of the two frames and $(\dv_u,\dv_v)$, $(\dv^\ast_u,\dv^\ast_v)$ their \emph{curvature} connectors.
Because the curvature tensors are symmetric, the curvature connectors must obey the identities
\begin{equation}\label{eq:connector_identities}
\dv_u\cdot\e_v=\dv_v\cdot\e_u\quad\text{and}\quad\dv_u^\ast\cdot\e_v^\ast=\dv_v^\ast\cdot\e_u^\ast.
\end{equation}

In \cite{sonnet:variational}, equations \eqref{eq:connectors} were employed to compute the tensor $\H$ as defined in \eqref{eq:H_definition} starting from the representation of $\R$ in \eqref{eq:rotationDyadic}. It was shown that the pure measures of drilling and bending envisioned in our theory take the form
\begin{equation}
\label{eq:pure_measures}
\wdr=4\lvert\V\cv^\ast-\R\cv\rvert^2\quad\text{and}\quad\wb=4(\lvert\V(\nablast\normal^{\ast})\rvert^2-\lvert\nablas\normal\rvert^2)^2.
\end{equation}

As we are here dealing with deformations between minimal surfaces,
we can give the bending energy a simpler form. A simple direct computation or an application of the Cayley-Hamilton theorem can be used to express the trace of the square of the
curvature tensor in terms of the mean curvature $H$ and the Gaussian curvature $K$ as
\begin{equation}
\lvert\nablas\normal\rvert^2
=\tr((\nablas\normal)^2)=2(2H^2-K).
\end{equation}
As for minimal surfaces $H=0$, we find by also using \eqref{eq:stretching} that
\begin{equation}\label{eq:wbCurvatures}
\wb=16\left(
\frac{\lvert\rv_{,u}^\ast\rvert^2}{\lvert\rv_{,u}\rvert^2}
K^\ast-K
\right)^2.
\end{equation}

\section{Neutral Modes of Minimal Surfaces}\label{sec:neutral}
We want to examine the three modes of deformation described in the preceding section within the scenario of Section~\ref{sec:kinematics}: two minimal surfaces are represented
by functions
\begin{equation}\label{eq:twoFs}
 F=\ex^{\Phi+i\chi}\quad\text{and}\quad F^\ast=\ex^{\Phi^\ast+i\chi^\ast},
\end{equation}
and
a deformation of one into the other is brought about by a conformal mapping
$h$ according to \eqref{eq:deformationMapping}.

In this case, the stretching measure \eqref{eq:w_s} can be found directly  from the expression \eqref{eq:stretchingU} for the stretching tensor $\U$,
which depends on $\Phi$, $\Phi^\ast$, and $h$. To find the drilling and bending measures given by \eqref{eq:pure_measures}, we need explicit expressions for the spin and curvature connectors on both $\surface$ and $\surface^\ast$.

We present in Appendix \ref{app:connectorsS} computations that deliver
in terms of the functions $\Phi$ and $\chi$
the connectors for $\surface$ used on the left of \eqref{eq:connectors}.

Regarding the connectors for $\surface^\ast$ used on the right of \eqref{eq:connectors}, the situation is slightly more involved because the
moving frame is parametrised in terms of the unstarred quantities $u$ and $v$. The relevant computations are presented in Appendix \ref{app:connectorsSstar},
and the resulting connectors turn out to depend on $\Phi^\ast$, $\chi^\ast$,
and $h$.

\subsection{Representation of the Deformation Measures}
\subsubsection{Stretching}
With \eqref{eq:stretchingU} we have that
\begin{equation}
\U-\proj=\left(\frac{\lvert\rv_{,u}^\ast\rvert}{\lvert\rv_{,u}\rvert}-1\right)\proj.
\end{equation}
As $\lvert \proj\rvert^2=\tr\left(\proj^2\right)=\tr\proj=2$, this means that the stretching measure \eqref{eq:w_s} takes the form
\begin{equation}\label{eq:w_s_min}
\ws=2\left(\frac{\lvert\rv_{,u}^\ast\rvert}{\lvert\rv_{,u}\rvert}-1\right)^2
=2\left(\ex^{\Phi^\ast-\Phi}\hpn\frac{1+\hn^2}{1+\wn^2}-1\right)^2,
\end{equation}
where the last equality is obtained by use of \eqref{eq:Weierstrass_metric_vector_u_v} and \eqref{eq:rStarunorm}.

\subsubsection{Drilling}
To compute the drilling energy density as given by \eqref{eq:pure_measures}, we
combine the rotation $\R$ and the stretching $\V$ as in \eqref{eq:rotationDyadic} and \eqref{eq:stretching} with the connectors $\cv$ and $\cv^\ast$ given in
\eqref{eq:Cconnector} and \eqref{eq:Cstarconnector} to find that
\begin{equation}\label{eq:w_d_min}
\wdr=
-K(1+\wn^2)^2
\Biggl
\{
\left[\alpha_{,u}-\left(\ln\frac{(1+\hn^2)\hpn}{1+\wn^2}\right)_{,v}\right]^2
+
\left[\alpha_{,v}+\left(\ln\frac{(1+\hn^2)\hpn}{1+\wn^2}\right)_{,u}
\right]^2
\Biggl\}.
\end{equation}
Here, we have set
\begin{equation}
\alpha=\chi^\ast-\chi
\end{equation}
and used that 
\begin{equation}
 K=-\frac{4}{\lvert\rv_{,u}\rvert^2(1+\wn^2)^2},
\end{equation}
see \eqref{eq:Weierstrass_metric_vector_u_v} and 
 \eqref{eq:Weierstrass_curvature}.

\subsubsection{Bending}
To find the bending energy density as given by \eqref{eq:pure_measures}, we
first note that 
\begin{equation}
\lvert\nablas\normal\rvert^2=\tr[(\nablas\normal)\trans(\nablas\normal)]=\lvert\dv_u\rvert^2+\lvert\dv_v\rvert^2,
\end{equation}
see \eqref{eq:connectors}. The connectors $\dv_u$ and $\dv_v$ are given in explicit form in \eqref{eq:DuConnector} and \eqref{eq:DvConnector}.

Similarly,
\begin{equation}
\lvert\nablast\normal^\ast\rvert^2
=\lvert\dv^\ast_u\rvert^2+\lvert\dv^\ast_v\rvert^2
=\lvert\dv^\ast_{u^\ast}\rvert^2+\lvert\dv^\ast_{v^\ast}\rvert^2,
\end{equation}
where $\dv^\ast_{u^\ast}$ and $\dv^\ast_{v^\ast}$ are defined in
\eqref{eq:DoubleStarConnectors} and given by \eqref{eq:DuConnector}
and \eqref{eq:DvConnector} with all relevant quantities starred.
Noting that
\begin{equation}
\V(\nablast\normal^{\ast})=\frac{\lvert\rv^\ast_{,u}\rvert}{\lvert\rv_{,u}\rvert}\nablast\normal^{\ast},
\end{equation}
we obtain after some simplifications that
\begin{equation}\label{eq:w_b_min}
\wb=16K^2\left[
\left(
\frac{\hpn(1+\wn^2)}{1+\hn^2}
\right)^2-1
\right]^2.
\end{equation}

\subsection{Mode Neutrality}
We say that a deformation is \emph{neutral} with respect to one of the three deformation modes if the corresponding energy density is
zero under that deformation. We are particularly interested in deformations that
are \emph{globally neutral} and are now in a position to derive criteria for global neutrality with respect to the three modes.

\subsubsection{Bending}
The bending energy density \eqref{eq:w_b_min} is
zero if and only if
\begin{equation}
\left(
\frac{\hpn(1+\wn^2)}{1+\hn^2}
\right)^2=1.
\end{equation}
As all three terms in the bracket on the left-hand side are real and strictly positive, this is equivalent to
\begin{equation}\label{eq:neutral_h}
\hpn=\frac{1+\hn^2}{1+\wn^2},
\end{equation}
which is a non-linear complex differential equation for $h$. For reasons that will soon become apparent, we shall
simply call \emph{neutral} any solution of \eqref{eq:neutral_h}.
We defer a discussion of the general solution of
\eqref{eq:neutral_h} to Section~\ref{sec:moebius}, where \eqref{eq:neutral_h} will also be given a geometric interpretation as an \emph{area preserving} mapping. For now, we note that there are indeed solutions of \eqref{eq:neutral_h}, one of them being the identity $h(w)=w$, the case treated in \cite{sonnet:bending-neutral}.

\subsubsection{Drilling}
The drilling energy density \eqref{eq:w_d_min} is
zero if and only if
\begin{equation}\label{eq:w_d_zero}
\alpha_{,u}=\left(\ln\frac{(1+\hn^2)\hpn}{1+\wn^2}\right)_{,v}
\quad\text{and}\quad
\alpha_{,v}=-\left(\ln\frac{(1+\hn^2)\hpn}{1+\wn^2}\right)_{,u}.
\end{equation}
By taking partial derivatives, it can be seen that this implies
that
\begin{equation}\label{eq:drillingDeltas}
\Delta\alpha=0
\quad\text{and}\quad
\Delta\ln\frac{(1+\hn^2)\hpn}{1+\wn^2}=0.
\end{equation}
The first of these equations is automatically satisfied: harmonic functions remain so under conformal mappings~\cite[\S 52]{weinberger:partial}, and so
as $\Delta^\ast\chi^\ast=0$ also $\Delta\chi^\ast=0$ and thus
$\Delta\alpha=\Delta(\chi^\ast-\chi)=0$.\footnote{It is perhaps worth recalling that here $\Delta^\ast$ denotes the Laplacian with respect to the starred variables $(u^\ast,v^\ast)$, whereas $\Delta$ is the Laplacian in the $(u,v)$ variables. To apply $\Delta$ to a starred function, such as $\chi^\ast$ or $\alpha^\ast$, the latter must be regarded as functions of $(u,v)$, upon use of \eqref{eq:u_v_star}.}

Regarding the second equation, note that for an analytic
function $f(z)=\lvert f(z)\rvert\ex^{i\arg f(z)}$ we have that
\begin{equation}\label{eq:lnf}
\ln f(z)=\ln\lvert f(z)\rvert + i\arg f(z).
\end{equation}
Thus wherever the logarithm of an analytic function is defined,
the logarithm of its modulus is the real part of an analytic function
and thus harmonic. Applied to $h'$, this means that
\begin{equation}
\Delta\ln\hpn=0,
\end{equation}
and thus a vanishing drilling energy density implies that
\begin{equation}
\Delta\ln \frac{1+\hn^2}{1+\wn^2}=0.
\end{equation}
An elementary computation shows that
\begin{equation}
\Delta \ln(1+\wn^2)=\frac{4}{(1+\wn^2)^2}
\quad\text{and so}\quad
\Delta^\ast\ln (1+\hn^2)=\frac{4}{(1+\hn^2)^2}.
\end{equation}
As $\Delta=\hpn^2\Delta^\ast$, see for example \cite[p. 548]{kythe:conformal}, we conclude that
\begin{equation}
\Delta\ln \frac{1+\hn^2}{1+\wn^2}=4\left(
\frac{\hpn^2}{(1+\hn^2)^2}-\frac{1}{(1+\wn^2)^2}
\right),
\end{equation}
which is zero if and only if $h$ satisfies \eqref{eq:neutral_h}, \emph{i.e.}, if $h$ is neutral. Thus a drilling-neutral deformation
is automatically also bending-neutral.

\subsubsection{Stretching}
The stretching energy density \eqref{eq:w_s_min} is
zero if and only if
\begin{equation}
\ex^{\Phi^\ast-\Phi}\hpn\frac{1+\hn^2}{1+\wn^2}=1.
\end{equation}
Dividing by the exponential and taking the logarithm shows
that this implies
\begin{equation}\label{eq:isometry_phi}
\ln\frac{(1+\hn^2)\hpn}{1+\wn^2}=\Phi-\Phi^\ast.
\end{equation}
With this, and setting
\begin{equation}\label{eq:phi}
 \varphi:=\Phi^\ast-\Phi,
\end{equation}
Equations \eqref{eq:w_d_zero} simply become
\begin{equation}
\alpha_{,u}=-\varphi_{,v}
\quad\text{and}\quad
\alpha_{,v}=\varphi_{,u}.
\end{equation}
That these equations are satisfied is a straightforward consequence
of the Cauchy-Rieman equations for $\ln F=\Phi+i\chi$,
$\ln F^\ast=\Phi^\ast+i\chi^\ast$, and $h$. Thus a stretching neutral deformation is automatically also drilling neutral (and hence bending neutral).

\subsubsection{Consequences of Mode Neutrality}
We have established a hierarchy of neutral modes, in particular we have that
\begin{equation}\label{eq:hierarchy}
\ws=0\quad\Rightarrow\quad
\wdr=0\quad\Rightarrow\quad
\wb=0.
\end{equation}
Furthermore, as $\wb=0$ if and only if Equation \eqref{eq:neutral_h} holds, whenever any mode is neutral, so is $h$. We note that
if there exist indeed isometric deformations of minimal surfaces
of the type \eqref{eq:deformationMapping}, then
\eqref{eq:hierarchy} shows that all modes are neutral and the
energy \eqref{eq:W} vanishes. Thus such isometries provide an
example of \emph{soft elasticity}.

In order to indentify some concrete neutral deformations, we now proceed to establish the restrictions placed by mode neutrality on
the minimal surfaces represented via $F$ and $F^\ast$. Naturally, the restrictions will be more severe when there are more neutral modes.

We start by considering $\ws=\wdr=\wb=0$. By using \eqref{eq:neutral_h} in \eqref{eq:isometry_phi}, we obtain that
with $\varphi$ as in \eqref{eq:phi} we have
\begin{equation}\label{eq:phi_hprime}
\varphi=-2\ln\hpn,\quad\text{i.e.,}\quad\ex^{\varphi}=\frac{1}{\hpn^2}.
\end{equation}
Furthermore, using \eqref{eq:neutral_h} in \eqref{eq:w_d_zero}
shows that
\begin{equation}
 \alpha_{,u}=(2\ln\hpn)_{,v}\quad\text{and}\quad
 \alpha_{,v}=-(2\ln\hpn)_{,u}.
\end{equation}
As $\ln\hpn$ and $\arg h'$ are conjugate harmonic functions, see \eqref{eq:lnf}, we conclude that
\begin{equation}\label{eq:alpha_arg}
 \alpha=\alpha_0-2\arg h'
\end{equation}
with some constant $\alpha_0$. As $F^\ast=\ex^{\varphi+i\alpha}F$,
see \eqref{eq:twoFs},
we find that
\begin{equation}\label{eq:isometryF}
F^\ast=\frac{1}{\hpn^2}\ex^{i(\alpha_0-2\arg h')}F
=\ex^{i\alpha_0}\frac{F}{h'^2}.
\end{equation}
When $h$ is the identity and so $h'=1$, this becomes
$F^\ast=\ex^{i\alpha_0}F$, which is the \emph{Bonnet transformation} of $F$, the same responsible for changing a catenoid into a helicoid (see, for example, pp.\,102 and 149 of \cite{dierkes:minimal}).

Let us now assume that $\ws\ne0$ but that still $\wdr=\wb=0$.
In that case, as $h$ is neutral, \eqref{eq:alpha_arg} still holds, but we need a replacement for \eqref{eq:phi_hprime}.  We
can use \eqref{eq:neutral_h} in \eqref{eq:w_s_min} to find
that
\begin{equation}\label{eq:w_s_bending_neutral}
\ws=
2\left(\ex^{\Phi^\ast-\Phi}\hpn^2-1\right)^2=
2(\lambda-1)^2,
\end{equation}
where we have set
\begin{equation}\label{eq:lambda}
 \lambda:=\hpn^2\ex^\varphi.
\end{equation}
Thus in this case the relation between $F$ and $F^\ast$
becomes
\begin{equation}\label{eq:drillingNeutralF}
F^\ast=\frac{\lambda}{\hpn^2}\ex^{i(\alpha_0-2\arg h')}F
=\lambda\ex^{i\alpha_0}\frac{F}{h'^2}.
\end{equation}

Finally, let us consider a situation in which $\ws\ne0\ne\wdr$ while still $\wb=0$. We can again use $\lambda$ as in \eqref{eq:lambda},
but as \eqref{eq:w_d_zero} no longer holds, we instead use
\eqref{eq:neutral_h} in the drilling energy density \eqref{eq:w_d_min} which thus becomes
\begin{align}
\wdr&=
-K(1+\wn^2)^2
\Biggl
\{
\left[\alpha_{,u}-\left(2\ln\hpn\right)_{,v}\right]^2
+
\left[\alpha_{,v}+\left(2\ln\hpn\right)_{,u}
\right]^2
\Biggl\}\nonumber\\
&=
-K(1+\wn^2)^2
\Biggl
\{
\left[\alpha_{,u}+\left(2\arg h'\right)_{,u}\right]^2
+
\left[\alpha_{,v}+\left(2\arg h'\right)_{,v}
\right]^2
\Biggl\}\nonumber\\
&=
-K(1+\wn^2)^2(\nabla\beta)^2,\label{eq:w_d_bending_neutral}
\end{align}
where we have once again used that $\ln\hpn$ and $\arg h'$ are conjugate harmonic functions and defined a harmonic function $\beta$
via
\begin{equation}
\beta:=\alpha+2\arg h'.
\end{equation}
As now $\alpha=\beta-2\arg h'$, we see that the function $\beta$ replaces the constant $\alpha_0$ in \eqref{eq:alpha_arg} and
the relation between $F$ and $F^\ast$
now is
\begin{equation}\label{eq:bendingNeutralF}
F^\ast=\lambda\ex^{i\beta}\frac{F}{h'^2}.
\end{equation}

To summarise, for a deformation to be \emph{bending neutral},
we need $F^\ast$ to be related to $F$ via \eqref{eq:bendingNeutralF}
with an $h$ that satisfies \eqref{eq:neutral_h}, a harmonic function $\beta$, and a stretching $\lambda$. If furthermore $\beta$ is constant, then \eqref{eq:drillingNeutralF} holds, and the deformation
is also \emph{drilling neutral}. In this case, if $\lambda=1$,
then \eqref{eq:isometryF} is satisfied, and the deformation is also an \emph{isometry}.

\subsection{M\"obius Transformation}\label{sec:moebius}
A special deformation of the type \eqref{eq:deformationMapping} is
the \emph{Goursat transformation}~\cite{goursat:transformation}. It is obtained by using for $h$ a M\"obius transformation\footnote{An effective visual representation of \eqref{eq:moebius} can be found in \cite{arnold:moebious}.}
\begin{equation}\label{eq:moebius}
h(w)=\frac{aw+b}{cw+d}\quad\text{with}\quad ad-bc\ne 0,
\quad a,b,c,d \in\mathbb{C},
\end{equation}
and by defining
\begin{equation}\label{eq:goursat}
F^\ast\left(h(w)\right):=\frac{F(w)}{(h'(w))^2},
\end{equation}
where, by \eqref{eq:moebius}, 
\begin{equation}\label{eq:goursat2}
h'=\frac{ad-bc}{(cw+d)^2}.
\end{equation}

It was already mentioned by Goursat that for the special
M\"obius transformations of the type
\begin{equation}\label{eq:moebiusSpecial}
h(w)=\frac{aw-\bar{c}}{cw+\bar{a}},
\end{equation}
which is obtained from \eqref{eq:moebius} by setting $b=-\bar{c}$ and $d=\bar{a}$, 
his transformation \eqref{eq:goursat} corresponds to a rigid rotation of the minimal surface \cite[p. 142]{goursat:transformation} (see also \S\,774 of \cite{nitsche:vorlesungen}).

The parameters featuring in \eqref{eq:moebius} or \eqref{eq:moebiusSpecial} are notoriously redundant; if all are multiplied by the same complex number, $h$ remains unaltered. Often  the following normalization conditions are adopted for  \eqref{eq:moebius} and \eqref{eq:moebiusSpecial},
\begin{equation}
\label{eq:moebius_normalizations}
ab-cd=1\quad\text{or}\quad |a|^2+|c|^2=1.
\end{equation}

It is a straightforward computation to show that $h$ in the form \eqref{eq:moebiusSpecial} satisfies
\eqref{eq:neutral_h}. Equation \eqref{eq:goursat} shows that
then \eqref{eq:isometryF} is satisfied with $\alpha_0=0$ and
so all the modes are neutral and the elastic energy \eqref{eq:W}
vanishes, as it should for a rigid rotation.

It can be shown that conformal mappings of the type
\eqref{eq:moebiusSpecial} are indeed the only solutions of
\eqref{eq:neutral_h}. To see this, we first need to provide a geometric interpretation of \eqref{eq:h_defintion}.

Letting $\sphere$ be the unit sphere in three-dimensional Euclidean space $\euclid$, we denote by $\pi$ the stereographic projection of $\sphere$ from its North pole $N$ onto the equatorial plane, which we identify with the complex plane $\mathbb{C}$. As shown, for example, in \cite[p.\,73]{needham:geometry}, the arc-length element $\dd s$ on $\sphere$ is related to the corresponding arc-length element $|\dd w|$ on $\mathbb{C}$ through the equation
\begin{equation}
	\label{eq:arc-length_element}
	\dd s=\frac{2}{1+|w|^2}|\dd w|.
\end{equation}
Rewriting \eqref{eq:neutral_h} in the equivalent form
\begin{equation}
	\label{eq:h_neutral_equivalent}
	\frac{2}{1+|w^\ast|^2}|\dd w^\ast|=\frac{2}{1+|w|^2}|\dd w|,
\end{equation}
where use has also been made of \eqref{eq:h_defintion}, we readily conclude that any solution $h$ of \eqref{eq:neutral_h} is indeed a locally \emph{area preserving} mapping of $\Omega$ onto $\Omega^\ast$.

Thus, a neutral $h$ as a mapping from $\mathbb{C}$ to $\mathbb{C}$
corresponds via the stereographic projection $\pi$ to
an area preserving mapping $\hat{h}$ of the Riemann sphere
$\hat{\mathbb{C}}$
onto itself.\footnote{The Riemann sphere $\hat{\mathbb{C}}$ is nothing but the complex plane $\mathbb{C}$ made compact by the addition of an extra point $\infty$, the point at \emph{infinity}. If, by the stereographic projection $\pi$, we  identify $\mathbb{C}$ with $\sphere\setminus\{N\}$, $\hat{\mathbb{C}}$ is then identified with the whole of $\sphere$ (see, for example, \cite[Sect.\,1.1]{jones:complex}).} As both $h$ and $\pi$ are
conformal, so is $\hat{h}$. But a classical theorem says that all conformal mappings $\hat{h}$ of $\hat{\mathbb{C}}$ into itself are represented by \eqref{eq:moebius} (see, for example \cite[p.\,269]{hilbert:geometry}) and it is a simple matter to show (see Appendix~\ref{sec:area}) that a M\"obius transformation as in \eqref{eq:moebius} is area preserving if and only if it can be set in the form \eqref{eq:moebiusSpecial}, and so it amounts to a rigid rotation of $\hat{\mathbb{C}}$. Said differently, any conformal, area preserving mapping of the unit sphere $\sphere$ onto itself
is an isometry generated by a uniform rotation in $\euclid$.\footnote{A somewhat weaker form of this result can also be given a direct proof, which is presented in \cite{sonnet:metric} among other metric implications for the kinematics of surfaces in three-dimensional space, much in the same spirit as \cite{seguin:coordinate}.}

In light of this, all modes of soft elasticity \eqref{eq:isometryF}
are either simply a rigid rotation of the minimal surface
when $\alpha_0=0$ or, for any other constant value of $\alpha_0$,
a rigid rotation combined with a Bonnet transformation.\footnote{This is just the same characterization for all isometries of minimal surfaces in three-dimensional space proved in \cite[p.\,159]{lawson:some}.} Thus, an isometry of a minimal surface is globally neutral, but it may fail to be a uniform rotation.

A further consequence of the fact that \eqref{eq:moebiusSpecial} is the only solution of \eqref{eq:neutral_h} is that, for bending neutral deformations, $h$ can be replaced by the identity $h(w)=w$, with no effective loss of generality, as done in \cite{sonnet:bending-neutral}; the only difference is that the \emph{spherical image} of $\surface^\ast$ would be the same as that of $\surface$, instead of being rotated in space.\footnote{By spherical image of a surface, we mean the image of the Gauss map $\normal:\surface\to\sphere$ associated with the surface.} Under a bending neutral deformation with $h$ as in  \eqref{eq:moebiusSpecial}, a minimal surface is changed into another with rotated spherical image, although one surface is \emph{not} necessarily the other rotated in space.

\subsection{Examples}\label{sec:examples}
Here, we shall illustrate with visual examples  deformations where different modes are selectively neutralized, including one where all modes are active. All images shown below depict surfaces represented in a Cartesian frame $\framee$ as in \eqref{eq:Weierstrass_representation} with $\Omega=\{w\in\mathbb{C}:1/e<|w|<e\}$ and $\Omega^\ast$ as in \eqref{eq:mappingh} for different Weierstrass functions $F$ and appropriate choices of $h$.

These are purely kinematic illustrations, where surfaces  are treated as geometric objects; as customary in the study of minimal surfaces, self intersections will be tolerated, although they would not be allowed for material surfaces.

Examples are listed according to a \emph{crescendo} of \emph{activated} modes.

\subsubsection{None}
\begin{figure}[] 
	\begin{subfigure}[b]{0.45\linewidth}
		\centering
		\includegraphics[width=\linewidth]{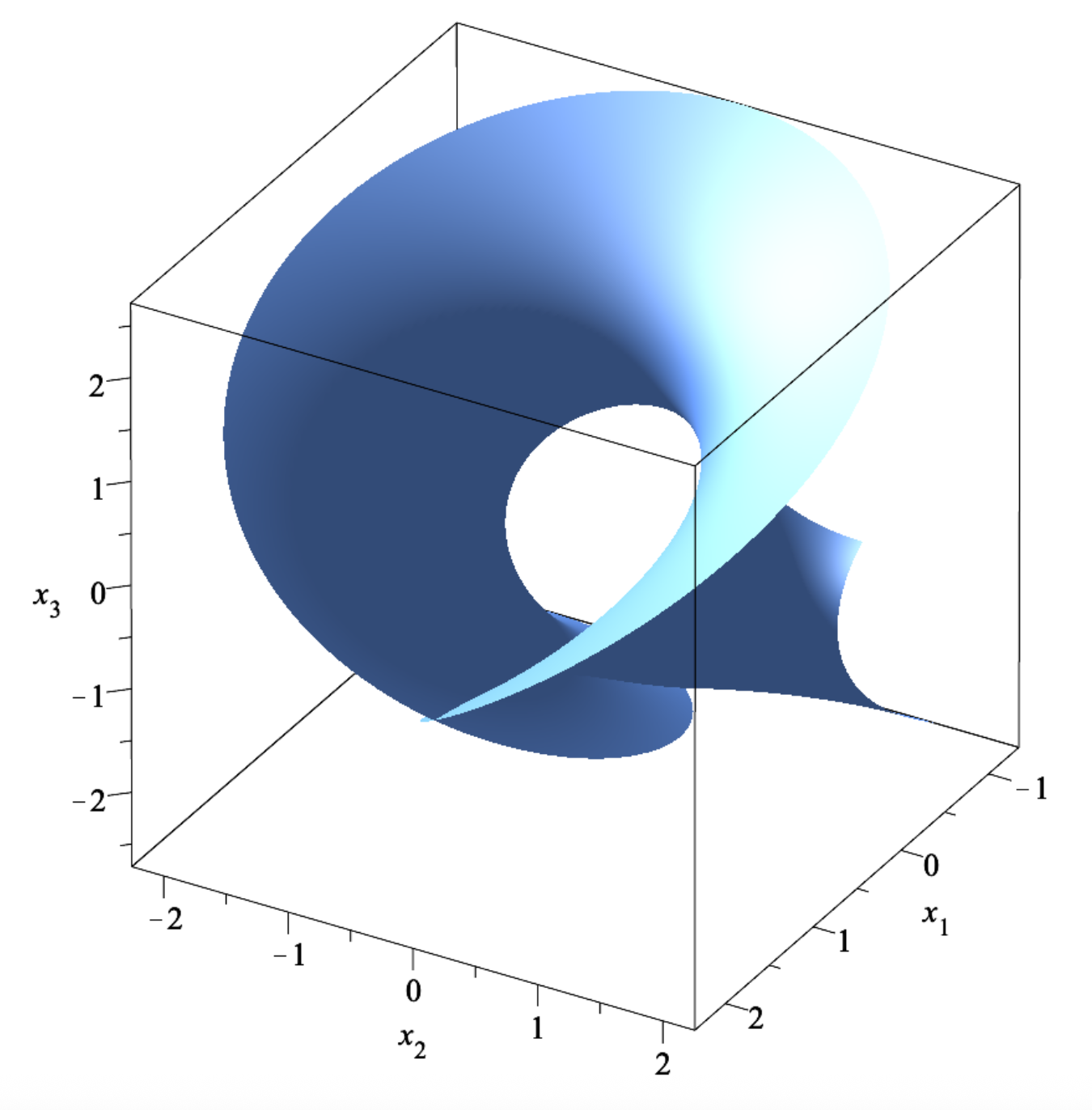}
		\caption{$F=\frac{1}{w}$}
		\label{fig:soft_elasticity_a}
	\end{subfigure}
	\quad
	\begin{subfigure}[b]{0.45\linewidth}
		\centering
		\includegraphics[width=\linewidth]{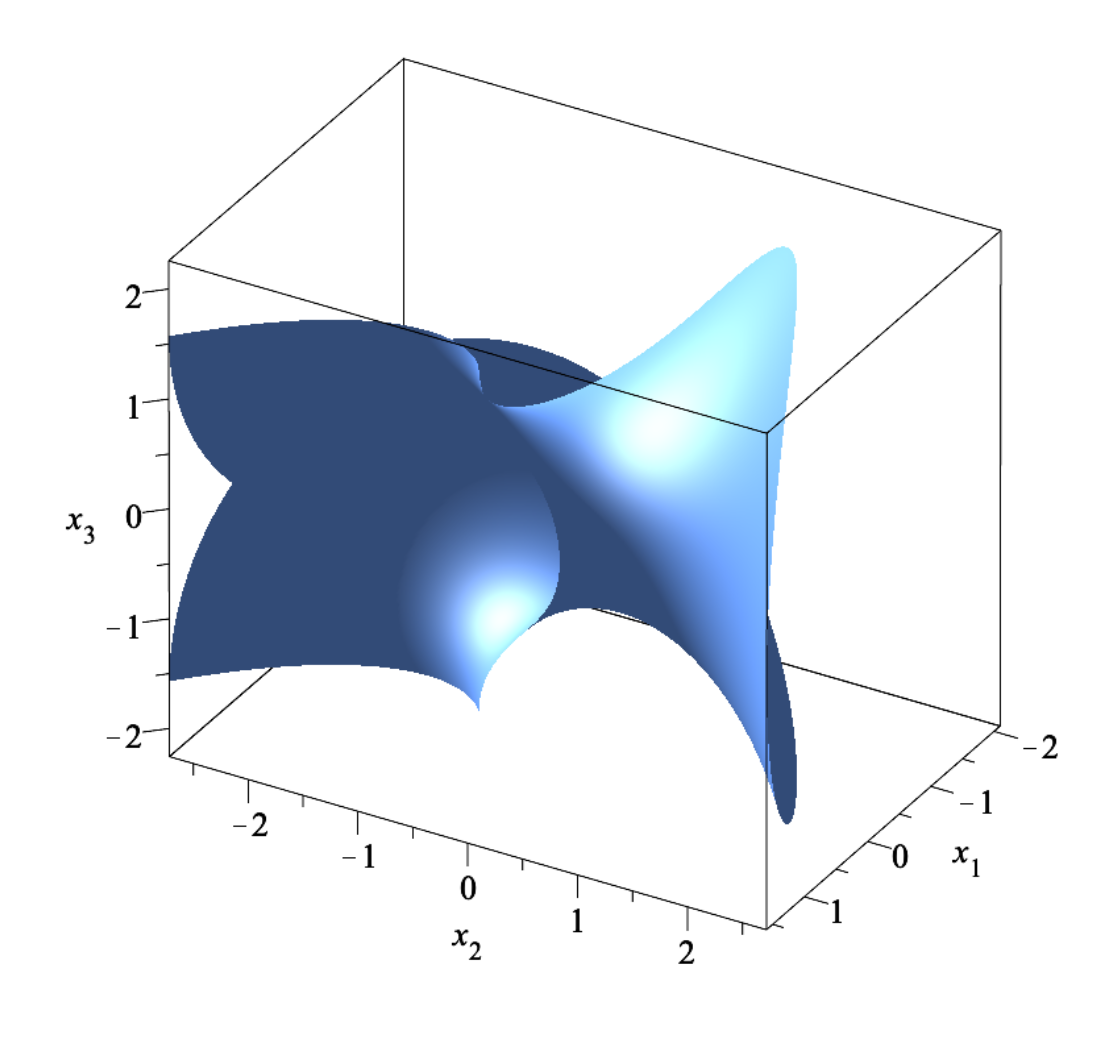}
		\caption{$\alpha_0=0$\quad$a=i$\quad$c=-1$}
		\label{fig:soft_elasticity_b}
	\end{subfigure}
	\quad
	\begin{subfigure}[b]{0.45\linewidth}
		\centering
		\includegraphics[width=\linewidth]{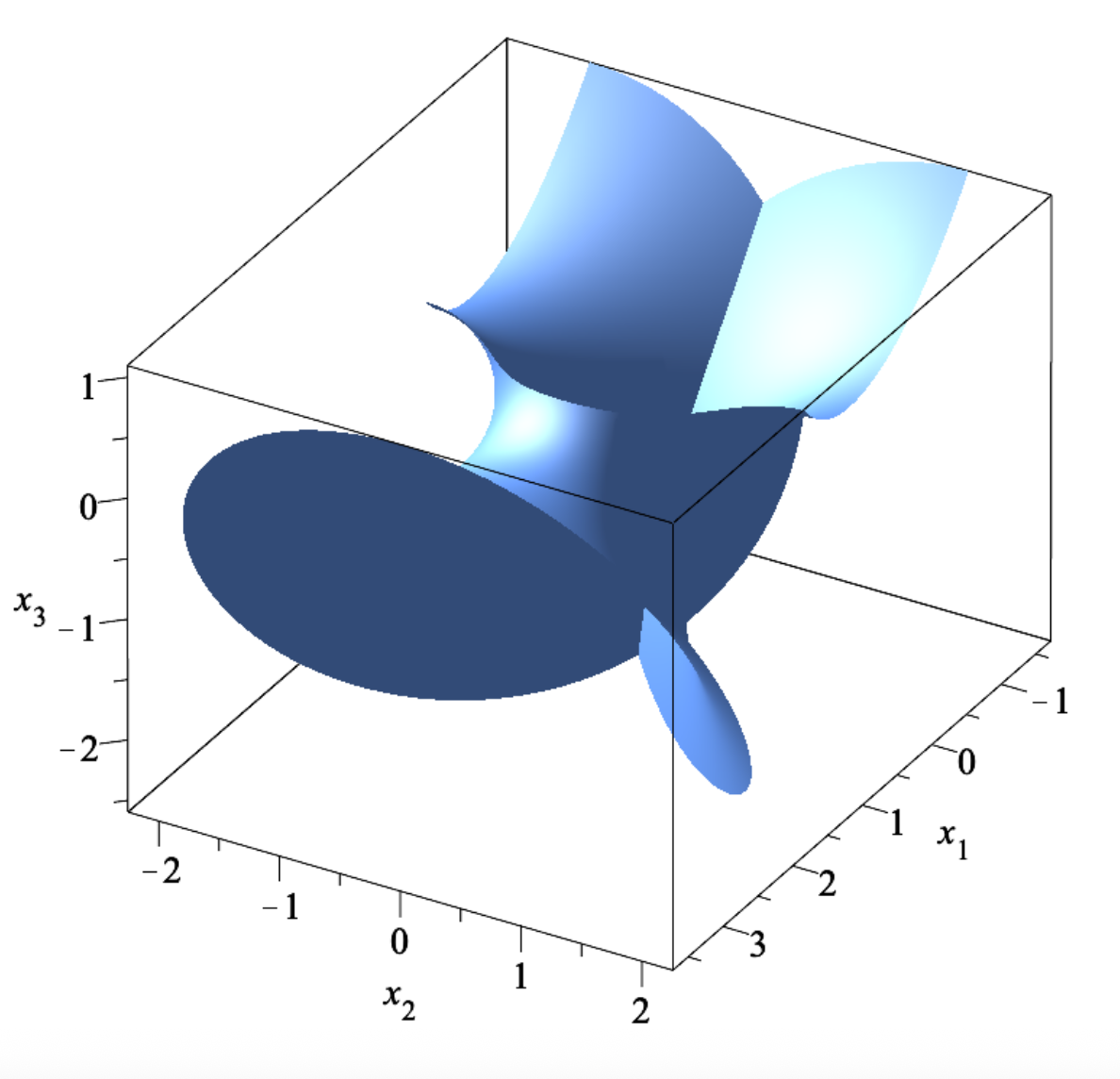}
		\caption{$\alpha_0=0$\quad$a=1$\quad$c=-1$}
		\label{fig:soft_elasticity_c}
	\end{subfigure}
		\quad
	\begin{subfigure}[b]{0.45\linewidth}
		\centering
		\includegraphics[width=\linewidth]{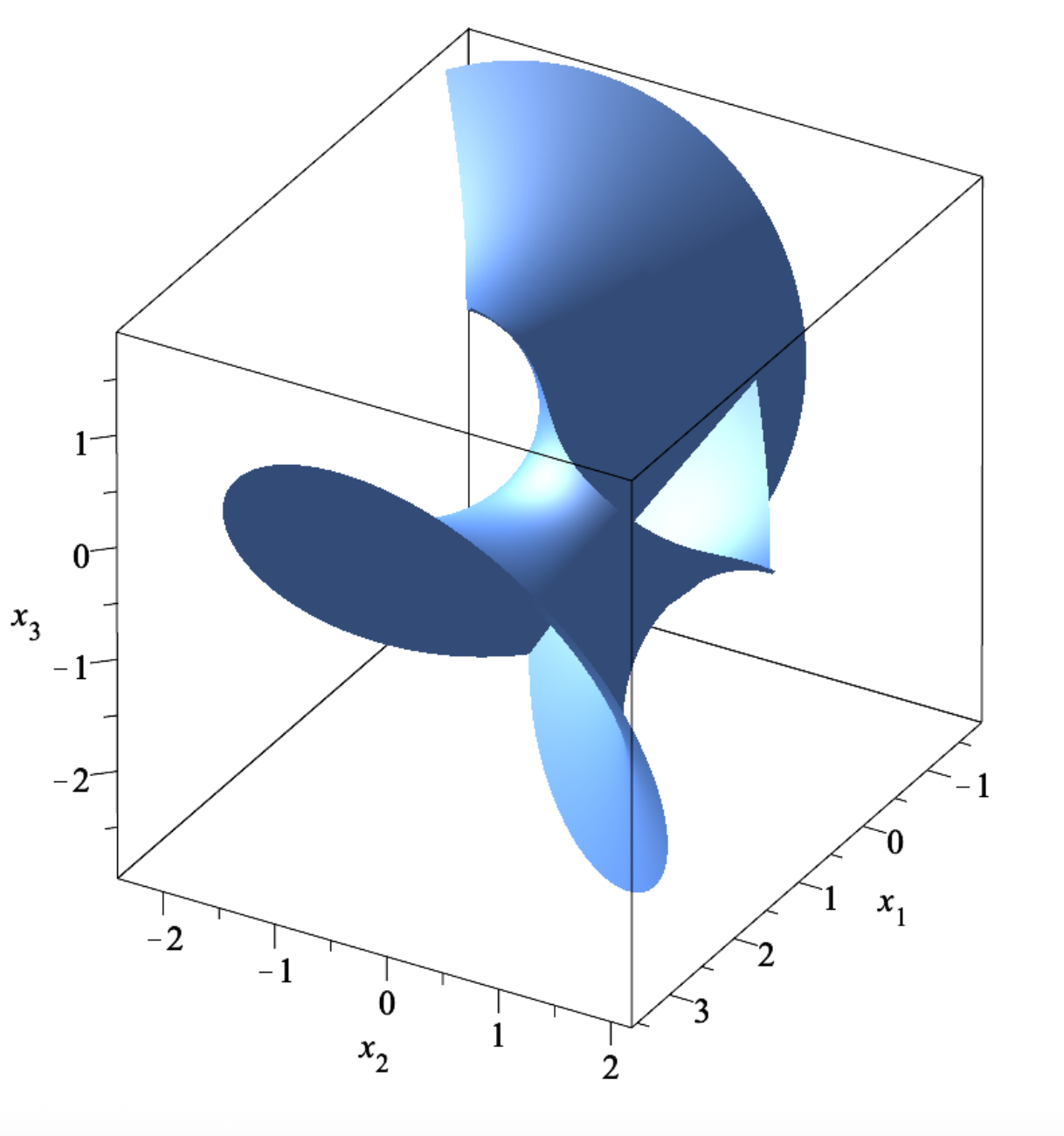}
		\caption{$\alpha_0=\frac{\pi}{6}$\quad$a=1$\quad$c=-1$}
		\label{fig:soft_elasticity_d}
	\end{subfigure}
	\quad
		\begin{subfigure}[b]{0.45\linewidth}
		\centering
		\includegraphics[width=\linewidth]{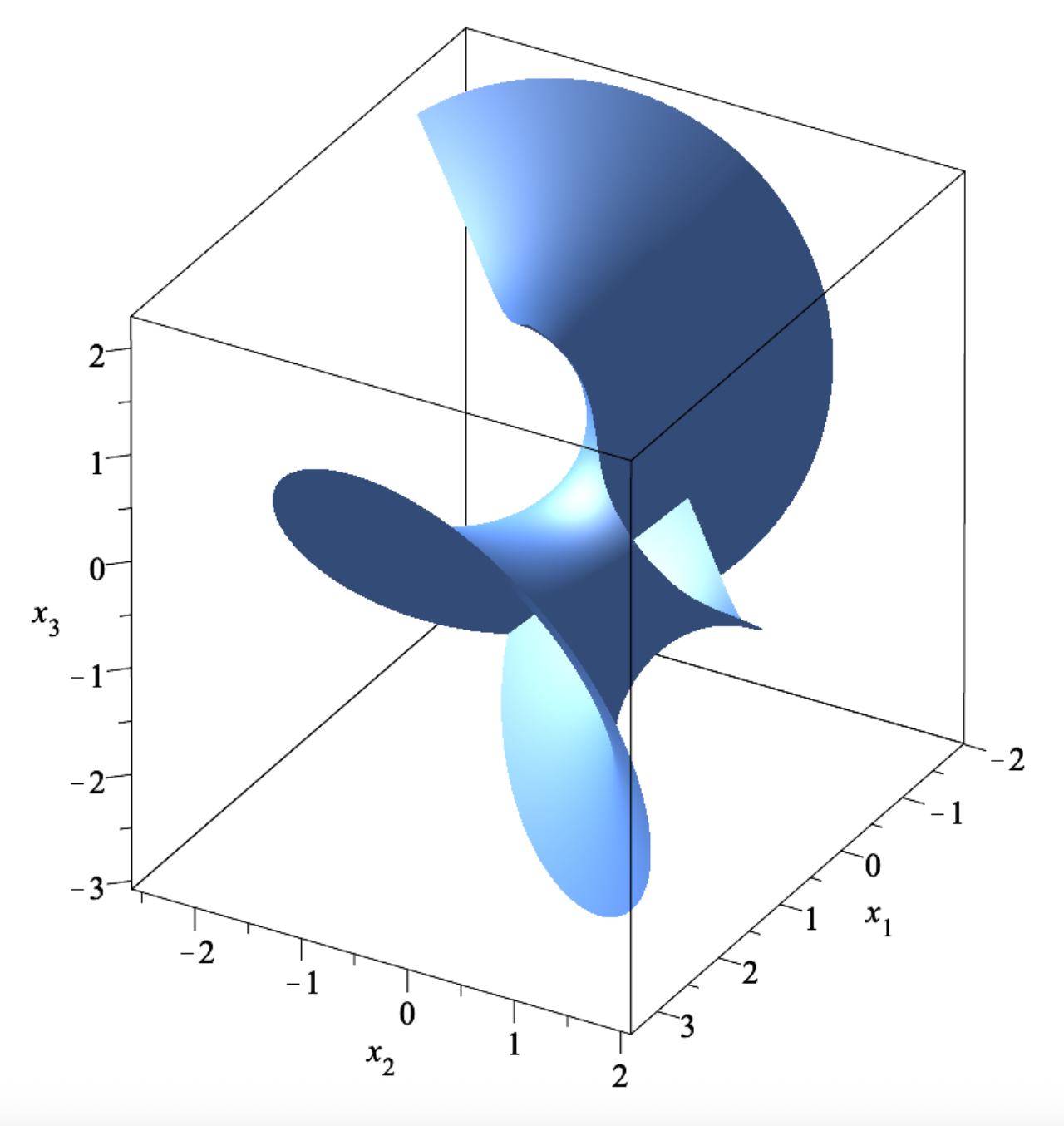}
		\caption{$\alpha_0=\frac{\pi}{4}$\quad$a=1$\quad$c=-1$}
		\label{fig:soft_elasticity_e}
	\end{subfigure}
	\quad
	\begin{subfigure}[b]{0.45\linewidth}
		\centering
		\includegraphics[width=\linewidth]{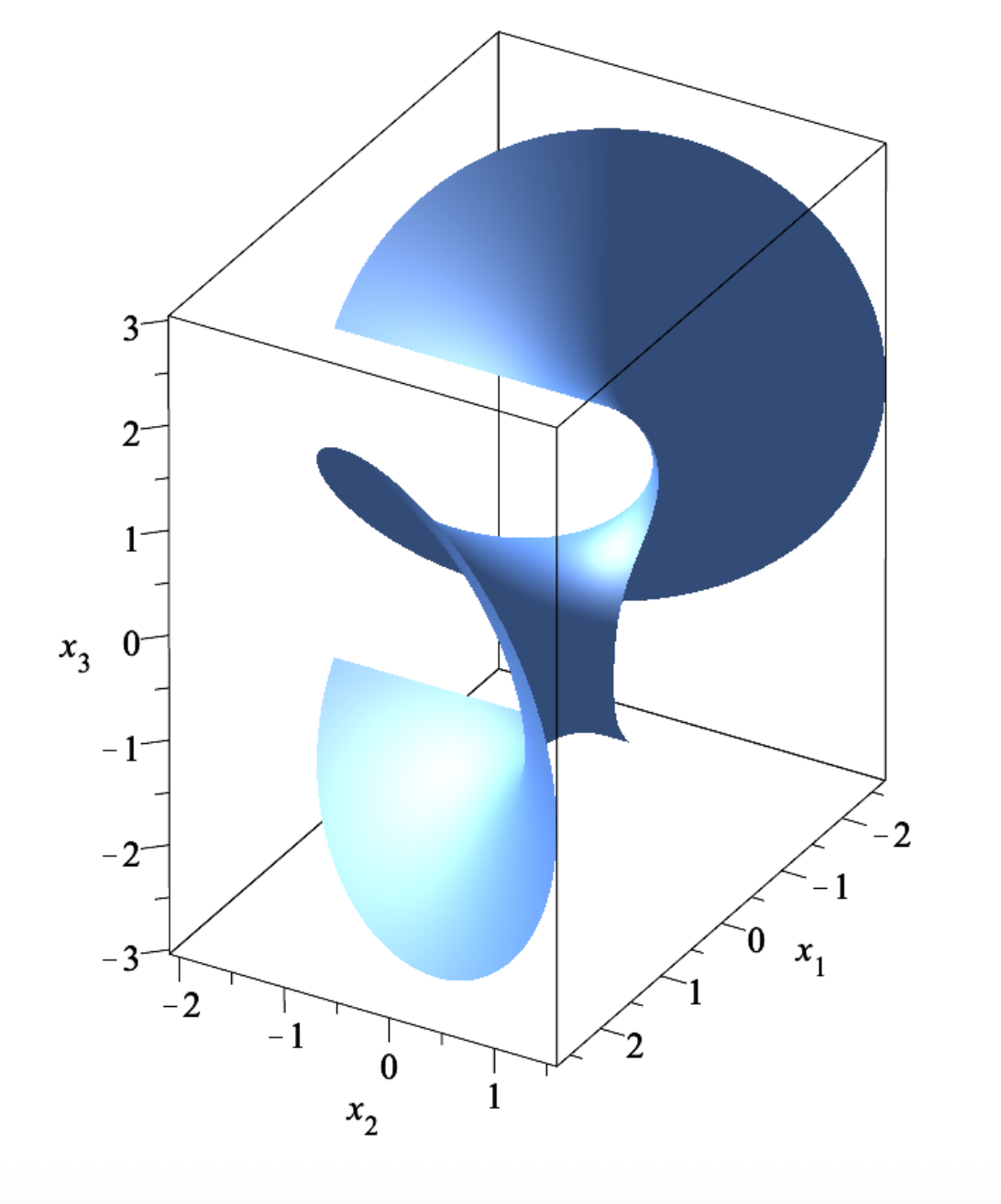}
		\caption{$\alpha_0=\frac{\pi}{2}$\quad$a=1$\quad$c=-1$}
		\label{fig:soft_elasticity_f}
	\end{subfigure}
	\caption{Bour's surface $\surface$ of index $m=1$ (represented by $F=1/w$) is deformed into a surface $\surface^\ast$ represented by $F^\ast$ in \eqref{eq:drillingNeutralF} with $\lambda=1$ and $h$ as in \eqref{eq:moebiusSpecial} for different values of $\alpha_0$ and complex parameters $a$, $c$. Panel (a) depicts the undeformed surface $\surface$; panels (b) and (c) depict rotations of $\surface$, whereas the surface $\surface^\ast$ in panels (d), (e), and (f) shows the compound effect of a rotation and a Bonnet transformation.}
	\label{fig:soft_elasticity}
\end{figure}

This is a case of soft elasticity, which we illustrate by setting $F=1/w$ and choosing $h$ as in \eqref{eq:moebiusSpecial}, so that $\surface$ is one of Bour's surfaces (with index $m=1$, as explained in \cite[p.\,156]{dierkes:minimal}). $F^\ast$ is given according to \eqref{eq:drillingNeutralF}, with $\lambda=1$ and different values of $\alpha_0$. Figure~\ref{fig:soft_elasticity} shows both $\surface$ and different instances of $\surface^\ast$, produced either by simple rotations of $\surface$ (when $\alpha_0=0$) or by the compound effect of a rotation and a Bonnet transformation (when $\alpha_0\neq0$).

Care should be used, both here and below, in representing $\surface^\ast$ thorough \eqref{eq:Weierstrass_representation}. By \eqref{eq:h_defintion} and \eqref{eq:drillingNeutralF}, we can write
\begin{subequations}
\begin{align}
	\rv^\ast(w^\ast)&=\Re\left(\frac12\int(1-w^{\ast2})F(w^\ast)\dd w^\ast\,\e_1+\frac{i}{2}\int(1+w^{\ast2})F(w^\ast)\dd w^\ast\,\e_2+\int w^\ast F(w^\ast)\dd w^\ast\,\e_3\right)\label{eq:Weierstrass_representation_starred_a}\\
	&=\Re\left(\frac12\int(1-h^{2})e^{i\alpha_0}\frac{F}{h'}\dd w\,\e_1+\frac{i}{2}\int(1+h^{2})e^{i\alpha_0}\frac{F}{h'}\dd w\,\e_2+\int w e^{i\alpha_0}\frac{F}{h'}\dd w\,\e_3\right).\label{eq:Weierstrass_representation_starred_b}
\end{align}
\end{subequations}
While the function in \eqref{eq:Weierstrass_representation_starred_a} is defined in $\Omega^\ast$, that in \eqref{eq:Weierstrass_representation_starred_b} is defined in $\Omega$, which in all applications in this section is an annulus centred at the origin of $\mathbb{C}$.

Needless to say, by \eqref{eq:hierarchy}, according to our theory,  all surfaces $\surface^\ast$ in Figs.~\ref{fig:soft_elasticity_b}-\ref{fig:soft_elasticity_f} are obtained from $\surface$ in Fig.~\ref{fig:soft_elasticity_a} at no elastic cost.

\subsubsection{Only stretching}\label{sec:stretching}
This case occurs when the deformation $\y$ is a uniform \emph{dilation} of 
$\surface$. Then $F^\ast=\lambda F$, with $\lambda>0$ and $h(w)=w$; by \eqref{eq:w_b_min} and \eqref{eq:w_d_bending_neutral}, both $\wb$ and $\wdr$ 
vanish, whereas $\ws$ does not if $\lambda\neq1$.

\subsubsection{Stretching and drilling}
This is the case studied in great detail in \cite{sonnet:bending-neutral}. $\surface$ is Enneper's surface represented by $F=1$ and $\surface^\ast$ is Bour's surface of index $m=3$, represented by $F^\ast=w^\ast$. Letting $h(w)=w$, from which it follows by \eqref{eq:w_b_min} that $\wb=0$, we can write $F^\ast$ as in \eqref{eq:bendingNeutralF} with $\lambda=|w|$ and $\beta=\arg w$. Simple calculations based on \eqref{eq:w_s_bending_neutral} and \eqref{eq:w_d_bending_neutral} show that 
\begin{equation}
	\label{eq:stretching_bending_energies}
	\ws=2(|w|-1)^2\quad\text{and}\quad\wdr=\frac{16}{(1+|w|^2)^2|w|^2},
\end{equation}
where \eqref{eq:Weierstrass_curvature} has also been used.

\begin{figure}[] 
	\begin{subfigure}[b]{0.45\linewidth}
		\centering
		\includegraphics[width=\linewidth]{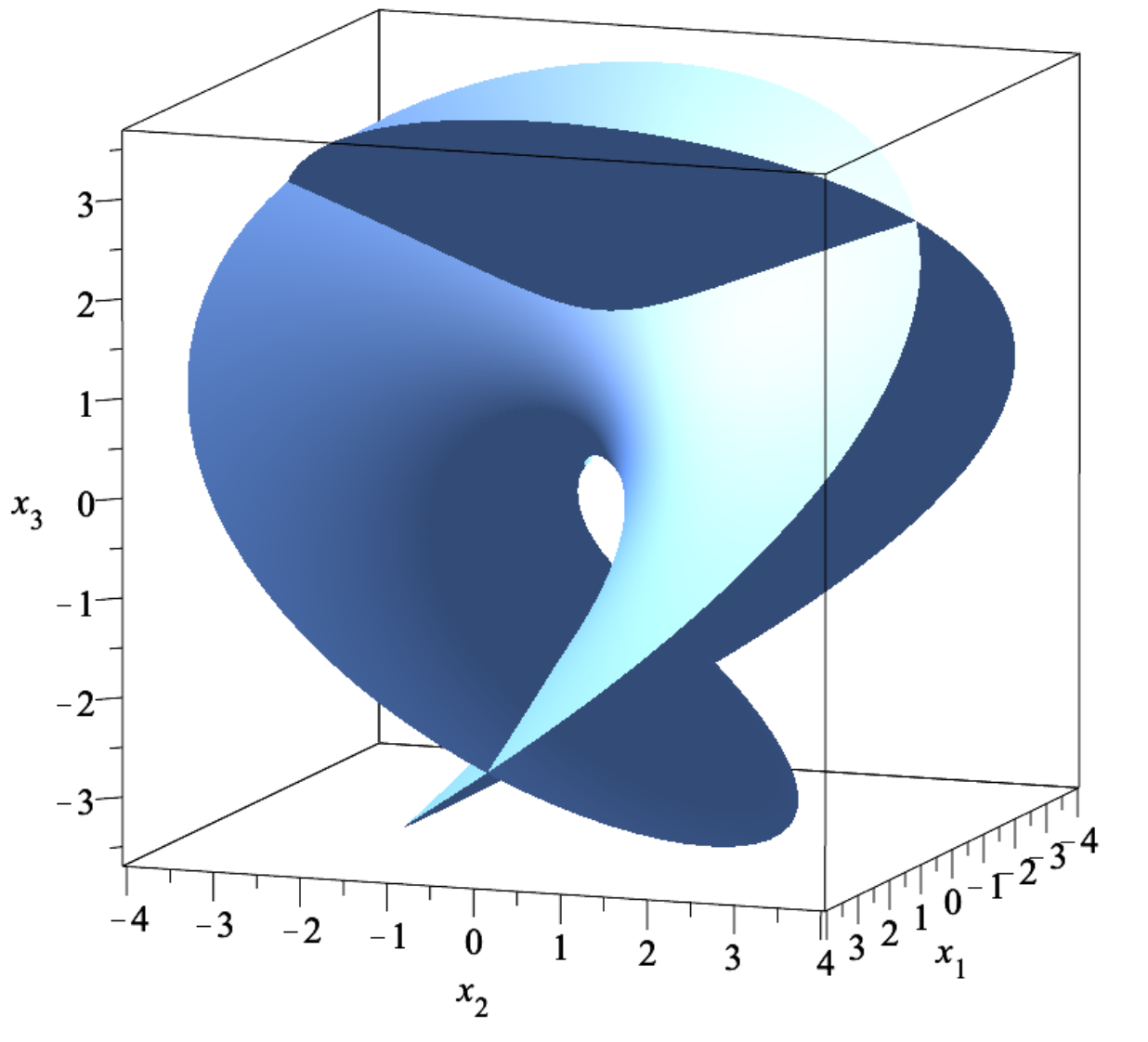}
		\caption{}
		\label{fig:bending_neutral_a}
	\end{subfigure}
	\quad
	\begin{subfigure}[b]{0.45\linewidth}
		\centering
		\includegraphics[width=\linewidth]{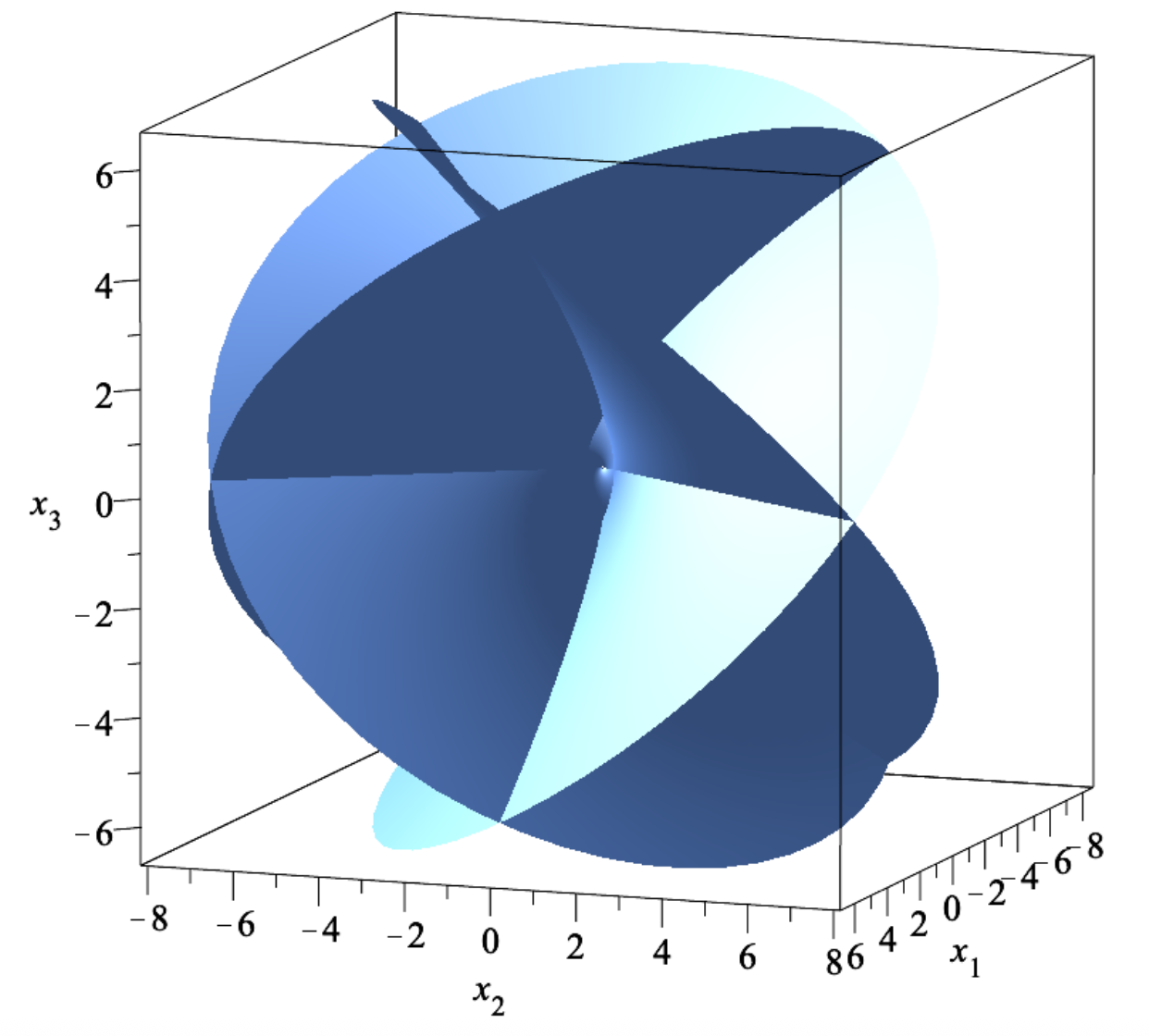}
		\caption{}
		\label{fig:bending_neutral_b}
	\end{subfigure}
	\caption{Enneper's surface (a) and Bour's surface of index $m=3$ (b). The latter is obtained by deforming the former with no bending energy; stretching and drilling energies are given by \eqref{eq:stretching_bending_energies}.}
	\label{fig:bending_neutral}
\end{figure}

For completeness. Fig.~\ref{fig:bending_neutral} shows both surfaces $\surface$ and $\surface^\ast$ corresponding to this case.

\subsubsection{All}\label{sec:all}
\begin{figure}[] 
	\begin{subfigure}[b]{0.45\linewidth}
		\centering
		\includegraphics[width=\linewidth]{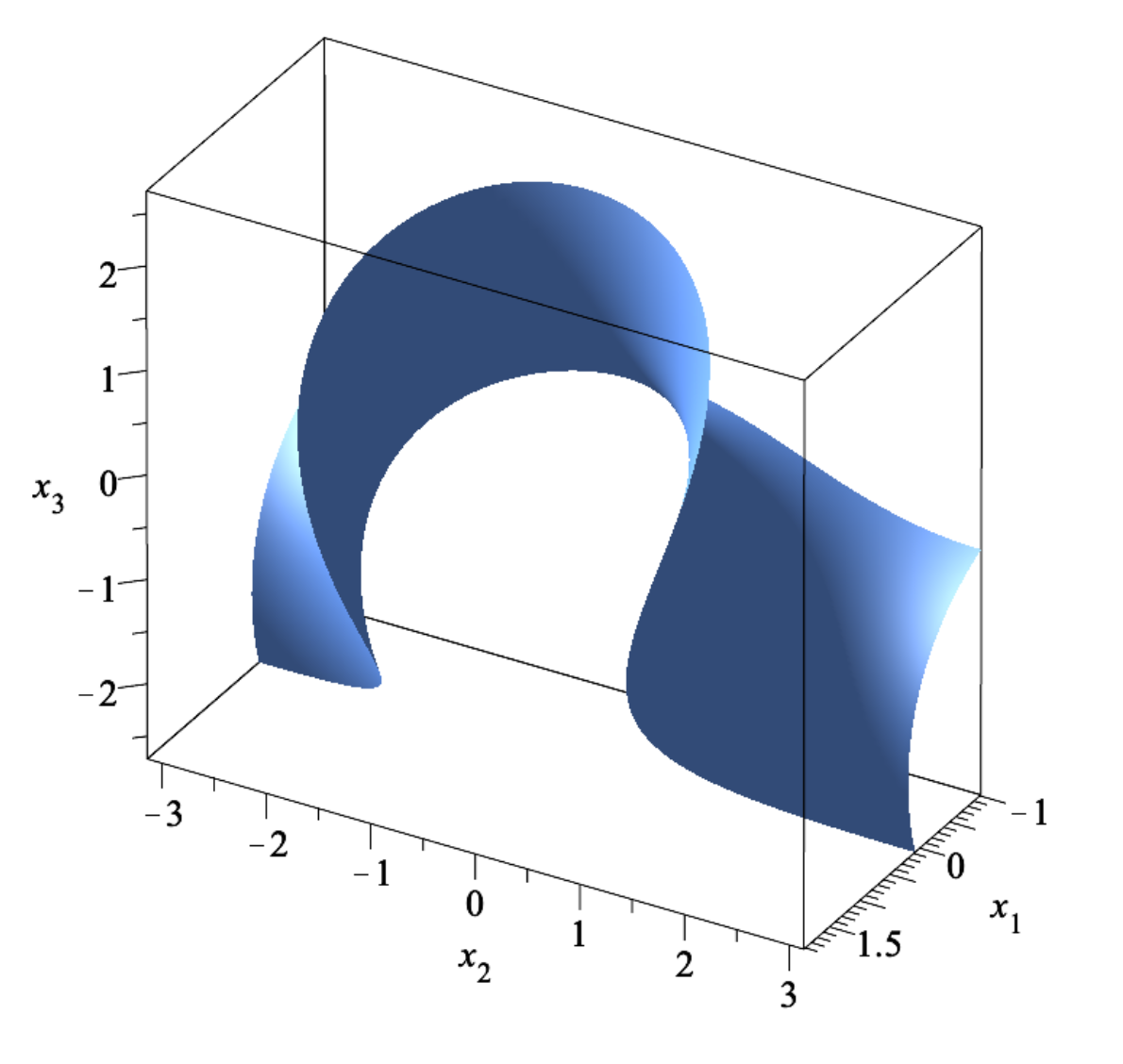}
		\caption{$\kappa=\frac12$}
		\label{fig:goursat_a}
	\end{subfigure}
	\quad
	\begin{subfigure}[b]{0.45\linewidth}
		\centering
		\includegraphics[width=\linewidth]{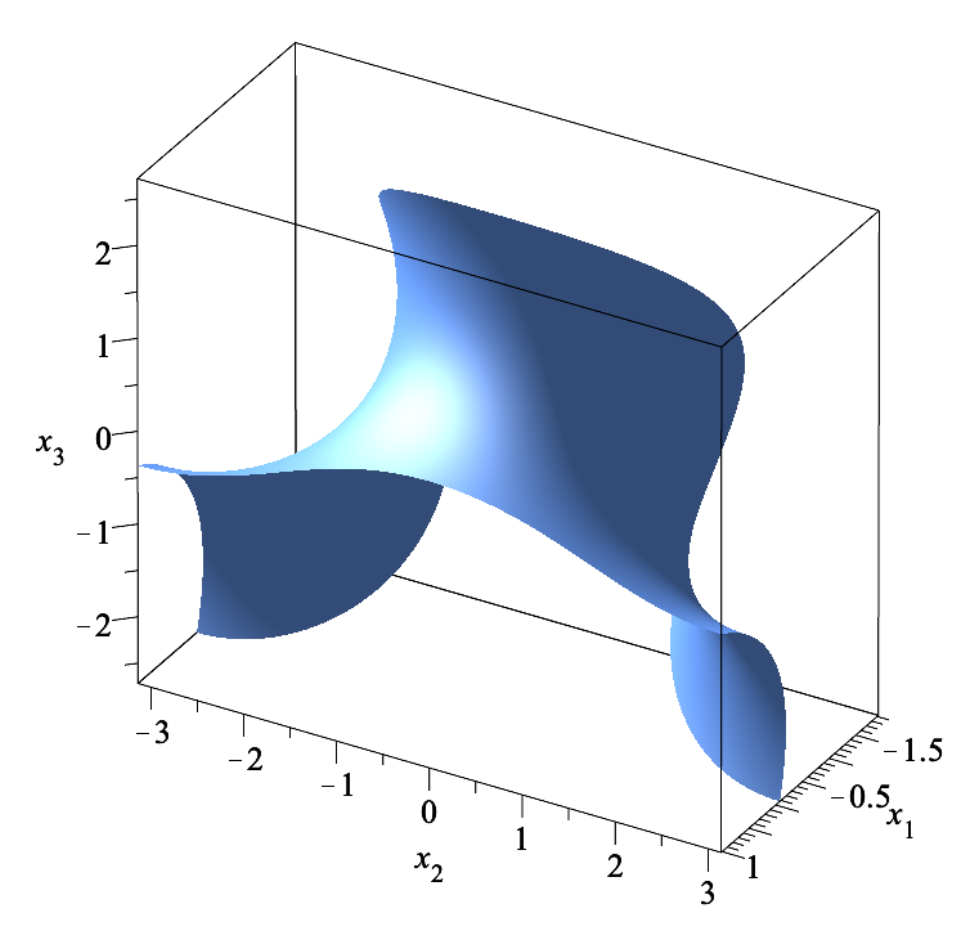}
		\caption{$\kappa=-\frac12$}
		\label{fig:goursat_b}
	\end{subfigure}
	\quad
		\begin{subfigure}[b]{0.45\linewidth}
		\centering
		\includegraphics[width=\linewidth]{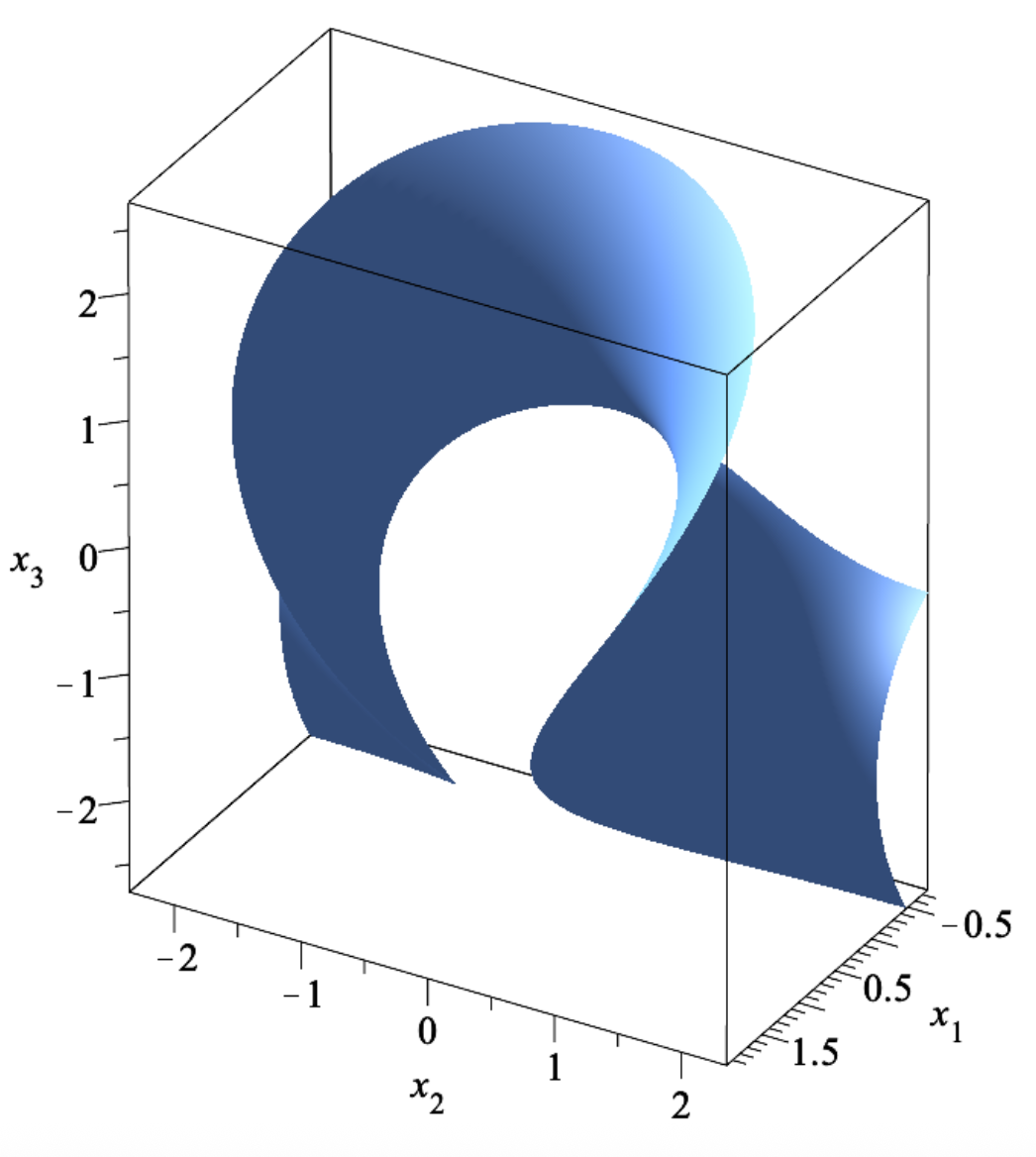}
		\caption{$\kappa=\frac23$}
		\label{fig:goursat_c}
	\end{subfigure}
	\quad
	\begin{subfigure}[b]{0.45\linewidth}
		\centering
		\includegraphics[width=\linewidth]{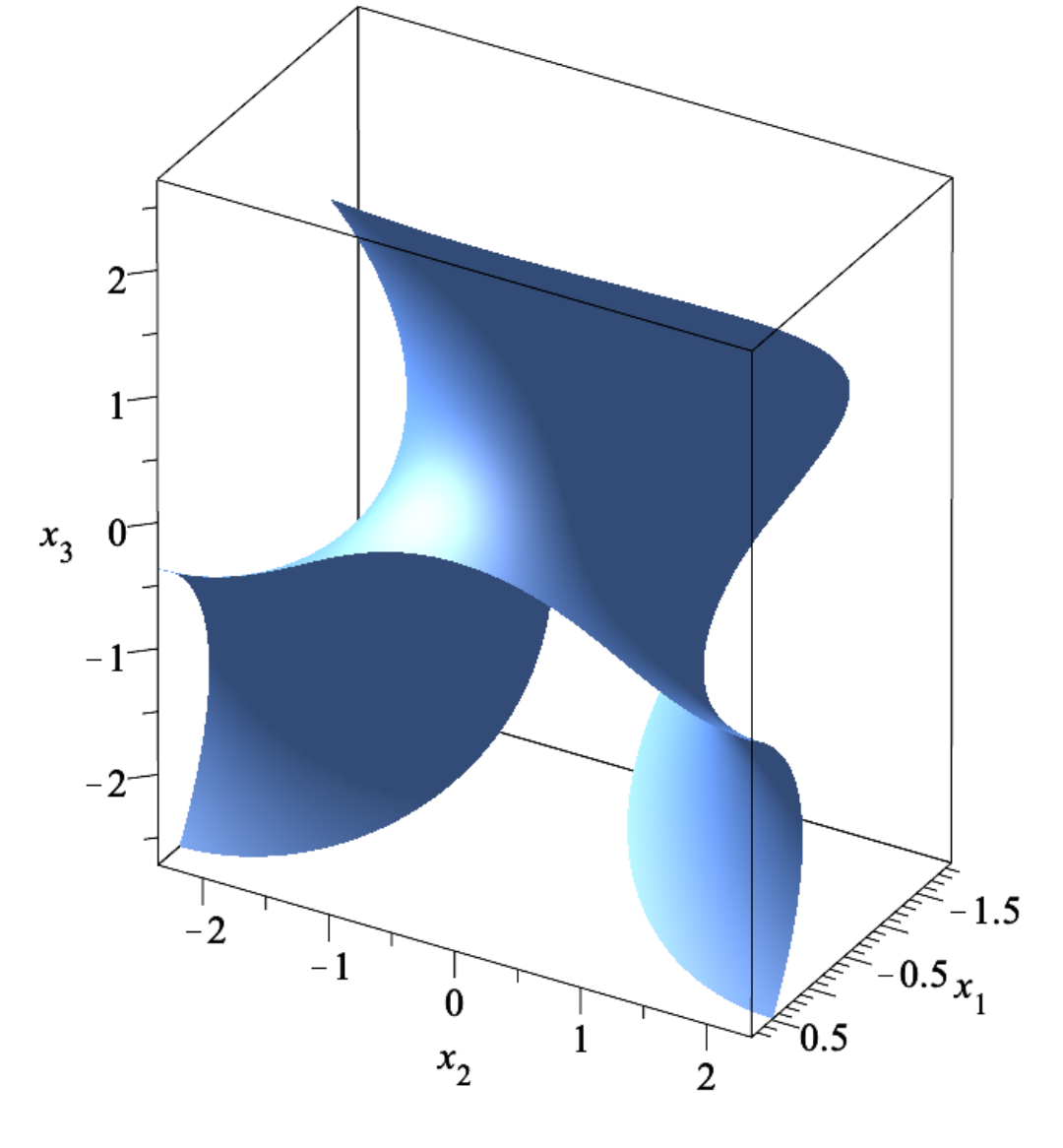}
		\caption{$\kappa=-\frac23$}
		\label{fig:goursat_d}
	\end{subfigure}
	\quad
		\begin{subfigure}[b]{0.45\linewidth}
		\centering
		\includegraphics[width=\linewidth]{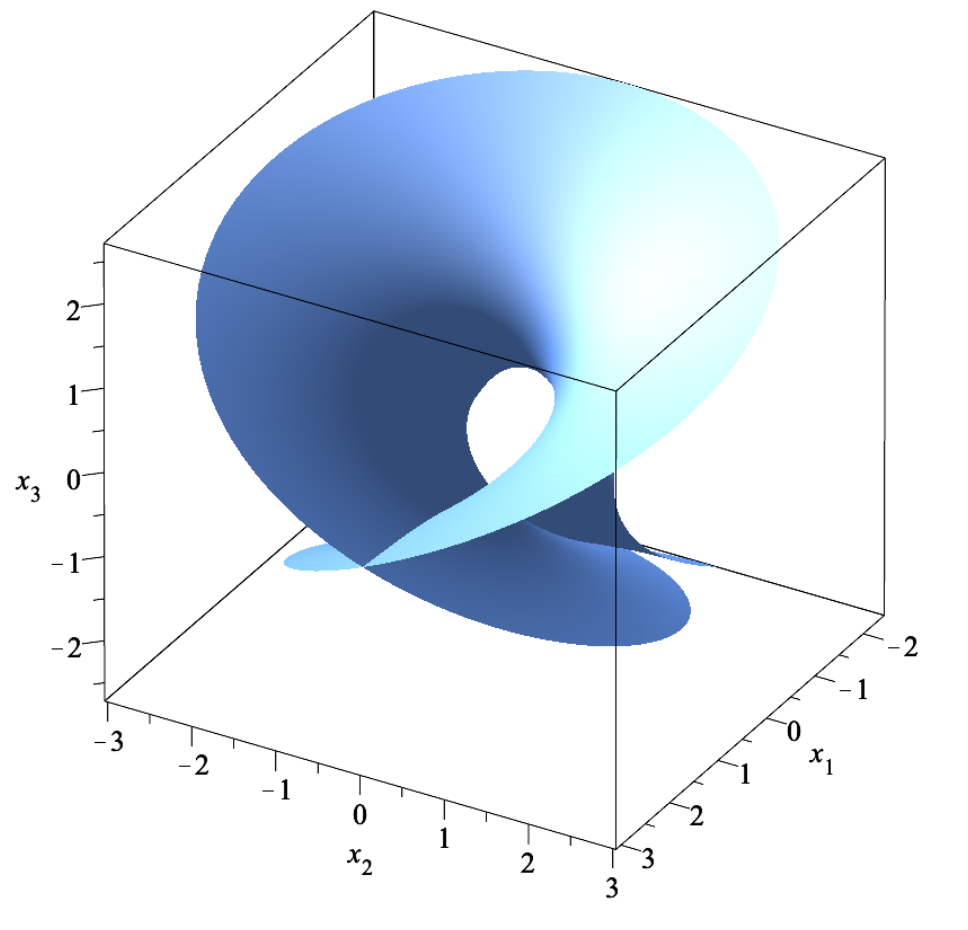}
		\caption{$\kappa=\frac32$}
		\label{fig:goursat_e}
	\end{subfigure}
	\quad\quad\quad
	\begin{subfigure}[b]{0.47\linewidth}
		\centering
		\includegraphics[width=\linewidth]{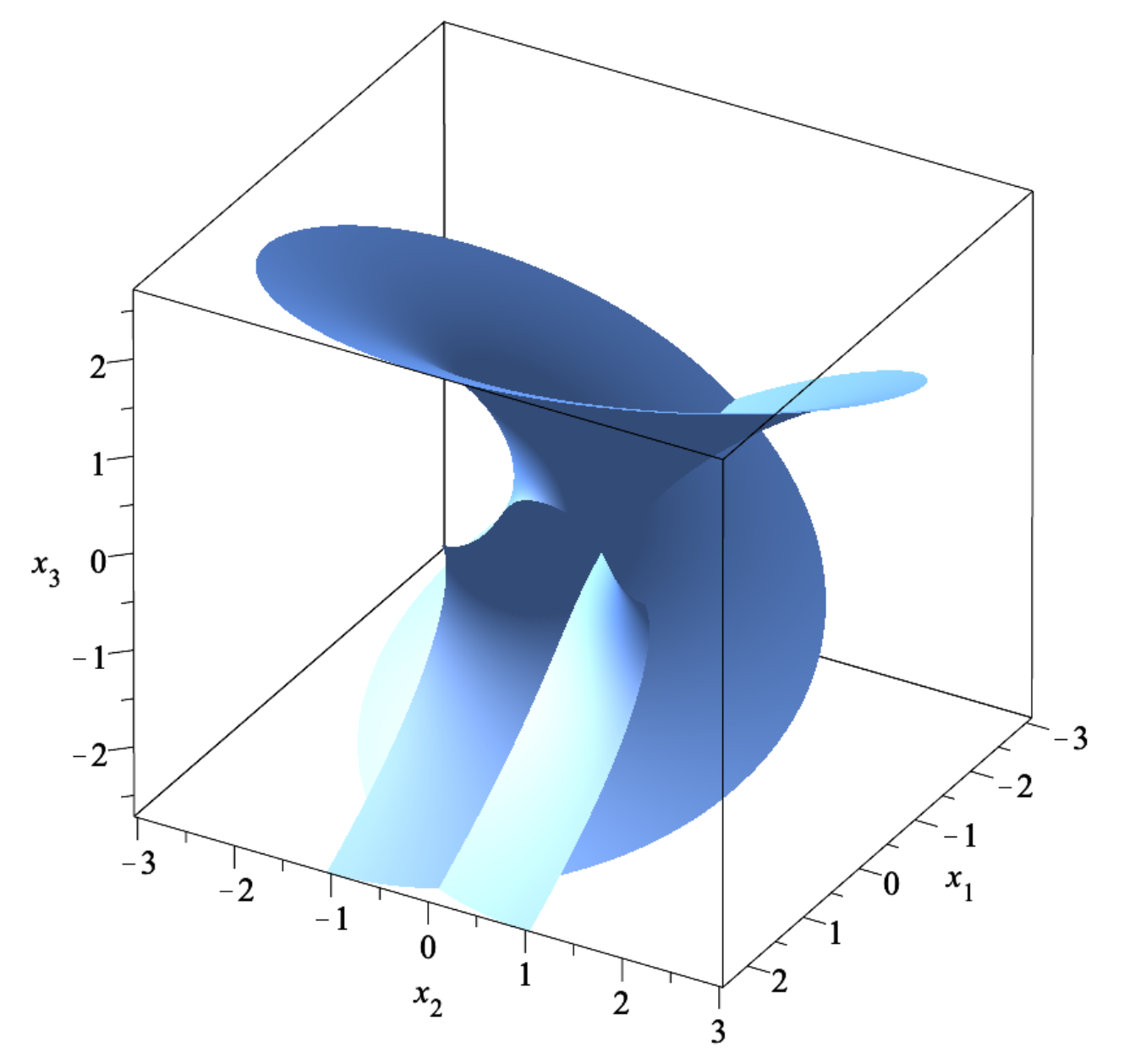}
		\caption{$\kappa=-\frac32$}
		\label{fig:goursat_f}
	\end{subfigure}
	\caption{Bour's surface of index $m=1$ depicted in Fig.~\ref{fig:soft_elasticity_a} is deformed via the Goursat transformation in \eqref{eq:goursat_kappa}, for different values of the parameter $\kappa$.}
	\label{fig:goursat}
\end{figure}

To see the theory applied in full, we consider deformations of the surface $\surface$ in Fig.~\ref{fig:soft_elasticity_a} (Bour's surface of index $m=1$) described by a special Goursat transformation, represented as in \eqref{eq:moebius} and \eqref{eq:goursat} with $a=\kappa\in\mathbb{R}$, $d=1$, and $b=c=0$ (see \cite[p.\,145]{goursat:transformation} and \cite[p.\,120]{dierkes:minimal}), that is,
\begin{equation}
	\label{eq:goursat_kappa}
	F^\ast(\kappa w)=\frac{F(w)}{\kappa^2},\quad h(w)=\kappa w.
\end{equation}
Figure~\ref{fig:goursat} shows the effect of this class of deformations for different values of the parameter $\kappa$ (both positive and negative). It is apparent that reverting the sign of $\kappa$ amounts to a rotation of $\pi$ about the axis $\e_3$.

Since $h$ in \eqref{eq:goursat_kappa} is \emph{not} neutral, none of the above simplified formulas for the energies of the independent modes apply here; we need to return to \eqref{eq:w_s_min}, \eqref{eq:w_d_min}, and \eqref{eq:w_b_min}. Simple calculations deliver the following expressions for the energies of the single modes,
\begin{subequations}
	\label{eq:goursat_kappa_energies}
	\begin{align}
	\ws&=2\left(1-\frac{1}{\kappa}\right)^2\left(\frac{1-\kappa|w|^2}{1+|w|^2}\right)^2,\label{eq:goursat_kappa_w_s}\\
	\wdr&=-4K\frac{(1-\kappa^2)^2|w|^2}{(1+\kappa^2|w|^2)^2},\label{eq:goursat_kappa_w_d}\\
	\wb&=16K^2\frac{(1-\kappa^2)^2(1-\kappa^2|w|^4)^2}{(1+\kappa^2|w|^2)^4},
	\end{align}
	where, by \eqref{eq:Weierstrass_curvature},
	\begin{equation}
		\label{eq:K_Bour_m=1}
		K=-\frac{16|w|^2}{(1+|w|^2)^4}.
	\end{equation}
\end{subequations}

\section{Conclusions}\label{sec:conclusion}
Stretching, drilling, and bending are the independent deformation modes occurring in a soft thin shell. For each, we proposed in \cite{sonnet:variational} a specific energy content, the combination of which provides the strain energy density of a second-grade variational theory.

Here, we considered all deformations of minimal surfaces into minimal surfaces and studied the deformation modes that selectively switch off one energy content; we called them \emph{neutral}.

A precise hierarchy has emerged. A stretching neutral mode, which is an \emph{isometry}, is also drilling neutral, and a drilling neutral mode is also bending neutral. Bending neutral modes, on the other hand, can be neither drilling  nor stretching neutral. The hardest, that is, most energetic deformations of a minimal surface into another harbour no neutral mode, whereas isometries are the softest of all, as they take place at no energy cost.

We hope that having characterized the energy modes of minimal surfaces may foster a better understanding of the energetics of soft thin shells. We found a means to classify all minimal surfaces relative to a reference one in terms of three different energies stored in them. It would be interesting to see whether this classification prompts any new pattern to emerge in the wide, fascinating realm of these surfaces. Our method actually goes beyond this realm: some of the conclusions reached here are indeed valid in a more general setting, as will be discussed in \cite{sonnet:metric}.


\appendix

\section{Connectors}\label{sec:appendix}

\subsection{Connectors on $\surface$}\label{app:connectorsS}
With the curvature tensor given in the form
\eqref{eq:Weierstrass_curvature_tensor}
we can easily read off the connectors $\dv_u$ and
$\dv_v$ in \eqref{eq:connectors}, we have
\begin{equation}\label{eq:DuConnector}
\dv_u=\frac{\ex^{\Phi}}{\lvert\rv_{,u}\rvert^2}(-\cos\chi\e_u+\sin\chi\e_v)
\end{equation}
and
\begin{equation}\label{eq:DvConnector}
\dv_v=\frac{\ex^{\Phi}}{\lvert\rv_{,u}\rvert^2}(\sin\chi\e_u+\cos\chi\e_v),
\end{equation}
where we have used that
$\displaystyle \frac{4\ex^{-\Phi}}{(1+\wn^2)^2}=\frac{\ex^{\Phi}}{\lvert\rv_{,u}\rvert^2}$.
To find $\cv$, note that
\begin{equation}
\cv\otimes\e_v=(\nablas\e_u)\trans-\dv_u\otimes\normal
\end{equation}
and so
\begin{equation}
\cv=(\nablas\e_u)\trans\e_v.
\end{equation}
Similar to \eqref{eq:Weierstrass_normal_dot} we have
\begin{equation}\label{eq:eu_dot}
\dot{\e}_u=\dot{u}\e_{u,u}+\dot{v}\e_{u,v}=(\nablas\e_u)\dot{\rv},
\end{equation}
and so we find with the help of \eqref{eq:Weierstrass_r_dot}
that
\begin{align}
\e_v\cdot\dot{\e}_u&=\dot{u}\e_v\cdot\e_{u,u}+\dot{v}\e_v\cdot\e_{u,v}=\e_v\cdot[(\nablas\e_u)\dot{\rv}]=
\dot{\rv}\cdot[(\nablas\e_u)\trans\e_v]=\dot{\rv}\cdot\cv\\
&=\lvert\rv_{,u}\rvert(\dot{u}\e_u\cdot\cv+\dot{v}\e_v\cdot\cv).
\end{align}
Therefore, because $\dot{u}$ and $\dot{v}$ are arbitrary,
\begin{equation}
 \cv=\frac{1}{\lvert\rv_{,u}\rvert}\left\{
(\e_v\cdot\e_{u,u})\e_u+(\e_v\cdot\e_{u,v})\e_v
\right\}.
\end{equation}
A straightforward if somewhat lengthy computation then yields
\begin{align}
\cv&=\frac{1}{\lvert\rv_{,u}\rvert}\left\{
\left(\chi_{,u}-\frac{2v}{1+\wn^2}\right)\e_u+
\left(\chi_{,v}+\frac{2u}{1+\wn^2}\right)\e_v
\right\}\\
&=\frac{1}{\lvert\rv_{,u}\rvert}\left\{
\left(\chi_{,u}-\left[\ln(1+\wn^2)\right]_{,v}\right)\e_u+
\left(\chi_{,v}+\left[\ln(1+\wn^2)\right]_{,u}\right)\e_v
\right\}.\label{eq:Cconnector}
\end{align}

\subsection{Connectors on $\surface^\ast$}\label{app:connectorsSstar}
The connectors $\dv^\ast_{u^\ast}$, $\dv^\ast_{v^\ast}$ and $\cv^\ast_\ast$ pertaining to $(\e^\ast_{u^\ast},\e^\ast_{v^\ast})$
satisfy
\begin{equation}\label{eq:DoubleStarConnectors}
  \begin{cases}
    \nablast\e_{u^\ast}^\ast=\e^\ast_{v^\ast}\otimes\cv^\ast_\ast+\normal^\ast\otimes\dv^\ast_{u^\ast},\\
    \nablast\e_{v^\ast}^\ast=-\e^\ast_{u^\ast}\otimes\cv^\ast_\ast+\normal^\ast\otimes\dv^\ast_{v^\ast},\\
    \nablast\normal^\ast=-\e^\ast_{u^\ast}\otimes\dv^\ast_{u^\ast}-\e^\ast_{v^\ast}\otimes\dv^\ast_{v^\ast}.
  \end{cases}
\end{equation}
They take exactly the same form as $\dv_u$, $\dv_v$, and $\cv$ as given
in \eqref{eq:DuConnector},  \eqref{eq:DvConnector}, and  \eqref{eq:Cconnector} with all relevant expressions starred. Note that while the two connectors $\cv^\ast$ and $\cv^\ast_\ast$ will usually differ from one another, the normal
$\normal^\ast$ and the curvature tensor $\nablast\normal^\ast$ are the same in \eqref{eq:connectors} and \eqref{eq:DoubleStarConnectors}.

We can use this last fact to determine $\dv^\ast_{u}$ and $\dv^\ast_{v}$. From \eqref{eq:connectors},
$\dv^\ast_{u}=-(\nablast\normal^\ast)\e^\ast_u$ and
$\dv^\ast_{v}=-(\nablast\normal^\ast)\e^\ast_v$, and so
by \eqref{eq:DoubleStarConnectors} we find that
\begin{align}
\dv^\ast_{u}&=(\e^\ast_{u^\ast}\cdot\e^\ast_{u})\,\dv^\ast_{u^\ast}+
(\e^\ast_{v^\ast}\cdot\e^\ast_{u})\,\dv^\ast_{v^\ast}\quad\text{and}\\
\dv^\ast_{v}&=(\e^\ast_{u^\ast}\cdot\e^\ast_{v})\,\dv^\ast_{u^\ast}+
(\e^\ast_{v^\ast}\cdot\e^\ast_{v})\,\dv^\ast_{v^\ast}.
\end{align}
Combining the starred versions of \eqref{eq:DuConnector} and \eqref{eq:DvConnector} with \eqref{eq:StarredBasesInverted} then yields that
\begin{align}
\dv^\ast_{u}=\frac{\ex^{\Phi^\ast}}{\lvert\rv^\ast_{,u}\rvert^2}
\{
&[(h_{u,v}^2-h_{u,u}^2)\cos\chi^\ast-2h_{u,u}h_{u,v}\sin\chi^\ast]\e^\ast_{u}\\\nonumber
+
&[2h_{u,u}h_{v,u}\cos\chi^\ast+(h_{u,u}^2-h_{u,v}^2)\sin\chi^\ast]\e^\ast_{v}
\}
\end{align}
and that
\begin{align}
\dv^\ast_{v}=\frac{\ex^{\Phi^\ast}}{\lvert\rv^\ast_{,u}\rvert^2}
\{
&[2h_{u,u}h_{v,u}\cos\chi^\ast+(h_{u,u}^2-h_{u,v}^2)\sin\chi^\ast]\e^\ast_{u}\\\nonumber
+
&[(h_{u,u}^2-h_{u,v}^2)\cos\chi^\ast+2h_{u,u}h_{u,v}\sin\chi^\ast]\e^\ast_{v}
\}.
\end{align}

We know that $\cv^\ast=-(\nablast\e^\ast_u)\trans\e^\ast_v$. To find the surface gradient of
$\e^\ast_u$ in the form \eqref{eq:estaru} note that%
\footnote{Use \eqref{eq:Weierstrass_r_dot} for $\rv^\ast(u,v)$,
$\dot{\rv}^\ast=\lvert\rv_{,u}^\ast\rvert(\dot{u}\e_u^\ast+\dot{v}\e_v^\ast)$, and
$\dot{f}=f_{,u}\dot{u}+f,_{,v}\dot{v}=(\nablast f)\dot{\rv}^\ast$.}
\begin{equation}\label{eq:nablasf}
 \nablast f=f_{,u}\frac{\e^\ast_u}{\lvert\rv^\ast_{,u}\rvert}+f_{,v}\frac{\e^\ast_v}{\lvert\rv^\ast_{,v}\rvert}
\end{equation}
for a scalar function $f$ defined on $\surface^\ast$ and parametrised in terms of $u$ and $v$. We
can thus write
\begin{equation}\label{eq:nablasEu}
\nablast\e^\ast_u=
\frac{h_{u,u}}{\hpn}\nablast\e^\ast_{u^\ast}+
\e^\ast_{u^\ast}\otimes\nablast\frac{h_{u,u}}{\hpn}+
\frac{h_{v,u}}{\hpn}\nablast\e^\ast_{v^\ast}+
\e^\ast_{v^\ast}\otimes\nablast\frac{h_{v,u}}{\hpn}.
\end{equation}
A direct computation using \eqref{eq:DoubleStarConnectors} shows that
\begin{align}
 \frac{h_{u,u}}{\hpn}(\nablast\e^\ast_{u^\ast})\trans\e^\ast_v+
 \frac{h_{v,u}}{\hpn}(\nablast\e^\ast_{v^\ast})\trans\e^\ast_v
 &=\frac{h_{u,u}}{\hpn}(\e^\ast_{v^\ast}\cdot\e^\ast_{v})\,\cv^\ast_\ast-
 \frac{h_{v,u}}{\hpn}(\e^\ast_{u^\ast}\cdot\e^\ast_{v})\,\cv^\ast_\ast\\
 &=\frac{h_{u,u}^2+h_{v,u}^2}{\hpn^2}\cv^\ast_\ast=\cv^\ast_\ast,
\end{align}
where \eqref{eq:StarredBasesInverted} has also been used.
Therefore from \eqref{eq:nablasEu} we obtain
\begin{align}
\cv^\ast=-(\nablast\e^\ast_u)\trans\e^\ast_v &=
\cv^\ast_\ast+(\e^\ast_{u^\ast}\cdot\e^\ast_{v})\,\nablast\frac{h_{u,u}}{\hpn}
+(\e^\ast_{v^\ast}\cdot\e^\ast_{v})\,\nablast\frac{h_{v,u}}{\hpn}\\
&=\cv^\ast_\ast-\frac{h_{v,u}}{\hpn}\nablast\frac{h_{u,u}}{\hpn}
+\frac{h_{u,u}}{\hpn}\nablast\frac{h_{v,u}}{\hpn},\label{eq:cstarOfcdoublestar}
\end{align}
again with the help of \eqref{eq:StarredBasesInverted}.

The starred form of \eqref{eq:Cconnector} gives $\cv^\ast_\ast$ in the basis
\eqref{eq:DoubleStarredBasis}. Once projected onto the basis
\eqref{eq:StarredBasis}, $\cv^\ast_\ast$ reads as
\begin{align}
\cv^\ast_\ast=\frac{1}{\lvert\rv^\ast_{,u^\ast}\rvert\hpn}
\Biggl
\{
&\left[h_{u,u}\chi^\ast_{,u^\ast}+h_{v,u}\chi^\ast_{,v^\ast}-2\frac{h_uh_{u,v}+h_vh_{v,v}}{1+\hn^2}\right]\e^\ast_{u}\nonumber\\
+
&\left[h_{u,v}\chi^\ast_{,u^\ast}+h_{v,v}\chi^\ast_{,v^\ast}+2\frac{h_uh_{u,u}+h_vh_{v,u}}{1+\hn^2}
\right]\e^\ast_{v}
\Biggl\}\\
=
\frac{1}{\lvert\rv^\ast_{,u}\rvert}
\Biggl
\{
&\left[\chi^\ast_{,u}-\left[\ln(1+\hn^2)\right]_{,v}\right]\e^\ast_{u}
+
\left[\chi^\ast_{,v}+\left[\ln(1+\hn^2)\right]_{,u}
\right]\e^\ast_{v}
\Biggl\},
\end{align}
where also the Cauchy-Riemann equations \eqref{eq:CRequations} have been used.

Then, using \eqref{eq:nablasf} to compute explicitly
the surface gradients in \eqref{eq:cstarOfcdoublestar}, we find
\begin{align}
-\frac{h_{v,u}}{\hpn}\nablast\frac{h_{u,u}}{\hpn}
+\frac{h_{u,u}}{\hpn}\nablast\frac{h_{v,u}}{\hpn}&=\frac{1}{\lvert\rv^\ast_{,u}\rvert\hpn^2}\{
(h_{u,uu}h_{u,v}-h_{u,uv}h_{u,u})
\e^\ast_{u}\nonumber\\
&\qquad\qquad\ \ +(h_{u,uv}h_{u,v}-h_{u,vv}h_{u,u})\e^\ast_{v}\}\\
&=\frac{1}{\lvert\rv^\ast_{,u}\rvert}\{
-\left[\ln\hpn\right]_{,v}\e^\ast_{u}
+
\left[\ln\hpn\right]_{,u}
\e^\ast_{v}\},
\end{align}
where again repeated use of \eqref{eq:CRequations} has been made.
Adding up, we obtain that
\begin{equation}\label{eq:Cstarconnector}
\cv^\ast
=\frac{1}{\lvert\rv^\ast_{,u}\rvert}
\Biggl
\{
\left[\chi^\ast_{,u}-\left[\ln(1+\hn^2)\hpn\right]_{,v}\right]\e^\ast_{u}
+
\left[\chi^\ast_{,v}+\left[\ln(1+\hn^2)\hpn\right]_{,u}
\right]\e^\ast_{v}
\Biggl\}.
\end{equation}

\section{Area preserving M\"obius transformations}\label{sec:area}

We are looking for M\"obius transformations
\begin{equation}\label{eq:moebiusNormalised}
f(z)=\frac{az+b}{cz+d},\quad ad-bc=1,\quad a,b,c,d\in\mathbb{C}
\end{equation}
that leave the area of any measurable subset of the Riemann sphere invariant, that is,
they have to satisfy
\begin{equation}\label{eq:areaPreservation}
\lvert f'(z)\rvert=\frac{1+\lvert f(z)\rvert^2}{1+\lvert z\rvert^2}.
\end{equation}
Explicitly, we have that
\begin{equation}
\lvert f'(z)\rvert=\frac{1}{(cz+d)(\bar{c}\bar{z}+\bar{d})}\quad\text{and}\quad
\lvert f(z)\rvert^2=\frac{(az+b)(\bar{a}\bar{z}+\bar{b})}{(cz+d)(\bar{c}\bar{z}+\bar{d})}.
\end{equation}

Multiplying \eqref{eq:areaPreservation} by $1+\lvert z\rvert^2$ and dividing by
$\lvert f'(z)\rvert$ delivers
\begin{align}
1+z\bar{z}&=(cz+d)(\bar{c}\bar{z}+\bar{d})+(az+b)(\bar{a}\bar{z}+\bar{b})\\
&=b\bar{b}+d\bar{d}+
(a\bar{a}+c\bar{c})z\bar{z}+
(a\bar{b}+c\bar{d})z+
(\bar{a}b+\bar{c}d)\bar{z}.
\end{align}
Comparing coefficients of $z$ and $\bar{z}$ on both sides of the equation shows that we need
\begin{subequations}
\begin{align}
b\bar{b}+d\bar{d}=&\ 1,\label{eq:area_a}\\
a\bar{a}+c\bar{c}=&\ 1,\label{eq:area_b}\\
a\bar{b}+c\bar{d}=&\ 0=\bar{a}b+\bar{c}d\label{eq:area_c}
\end{align}
together with the normalisation
\begin{equation}
 ad-bc=1.\label{eq:area_d}
\end{equation}
\end{subequations}
Subtracting from the sum of \eqref{eq:area_a} and \eqref{eq:area_b}
both \eqref{eq:area_d} and its complex conjugate delivers
\begin{align}
0&=a\bar{a}+b\bar{b}+c\bar{c}+d\bar{d}
-ad-\bar{a}\bar{d}+bc+\bar{b}\bar{c}\\
&=(a-\bar{d})(\bar{a}-d)+(b+\bar{c})(\bar{b}+c)
=\lvert a-\bar{d}\rvert^2+\lvert b+\bar{c}\rvert^2.
\end{align}
This shows that
\begin{equation}\label{eq:area_solution}
d=\bar{a}\quad\text{and}\quad b=-\bar{c},
\end{equation}
which also satisfies \eqref{eq:area_c}.

Using \eqref{eq:area_solution} in  \eqref{eq:moebiusNormalised} gives
\begin{equation}
f(z)=\frac{az-\bar{c}}{cz+\bar{a}},\quad \lvert a\rvert^2+\lvert c\rvert^2=1,\quad a,c\in\mathbb{C}.
\end{equation}
So the M\"obius transformations that leave the area invariant are
exactly those that deliver rigid rotations of the Riemann sphere.




%

\end{document}